\documentclass[a4paper,10pt,reqno]{amsart}

\usepackage{amssymb,amsmath,amscd,amsthm,hyperref}

\title[Two-term asymptotics for the kinetic energy operator]{Two-term spectral asymptotics for the Dirichlet pseudo-relativistic kinetic energy operator on a bounded domain}
\author{Sebastian Gottwald}
\date{}

\newtheorem{thm}{Theorem}[section]
 \newtheorem{cor}[thm]{Corollary}
 \newtheorem{lem}[thm]{Lemma}
 \newtheorem{prop}[thm]{Proposition}
 \theoremstyle{definition}
 
 \theoremstyle{remark}
 \newtheorem{rem}[thm]{Remark}
 
 \numberwithin{equation}{section}

\newcommand{\tr}{\textnormal{Tr}}
\newcommand{\supp}{\textnormal{supp\,}}

\newcommand{\h}{\hspace{1pt}}

\newcommand{\np}[1]{(#1)\hspace{-1pt}\underline{\hspace{5pt}}\h}
\newcommand{\npb}[1]{\big(#1\big)_{\hspace{-2pt}-}}
\newcommand{\npbb}[1]{\Big(#1\Big)_{\hspace{-2pt}-}}
\newcommand{\npzero}[1]{(#1)^{\hspace{0pt}\text{\tiny0}}\hspace{-5pt}
\underline{\hspace{5pt}}\h}

\newcommand{\defeq}{\mathrel{\mathop:}=}
\newcommand{\eqdef}{=\mathrel{\mathop:}}
\newcommand{\disp}[1]{\hspace{1pt}_{\text{\normalsize{#1}}}}

\usepackage[update,prepend]{epstopdf}

\begin{document}

\subjclass[2010]{35P20, 47G30}

\keywords{Spectral asymptotics, semiclassical analysis, Riesz means, Dirichlet condition, pseudo-relativistic operators}

\date{\today}

\begin{abstract}
Continuing the series of works following Weyl's one-term asymptotic formula for the counting function $N(\lambda)=\sum_{n=1}^\infty\npzero{\lambda_n{-}\lambda}$ of the eigenvalues of the Dirichlet Laplacian \cite{Weyl1912} and the much later found two-term expansion on domains with highly regular boundary by Ivri{\u\i} \cite{Ivrii1980} and Melrose \cite{Melrose1980}, we prove a two-term asymptotic expansion of the $N$-th Ces\`aro  mean of the eigenvalues of $A_m^\Omega\defeq\sqrt{-\Delta + m^2} - m$ for $m>0$ with Dirichlet boundary condition on a bounded domain $\Omega\subset\mathbb R^d$ for $d\geq 2$, extending a result by Frank and Geisinger \cite{Frank2013} for the fractional Laplacian ($m=0$) and improving upon the small-time asymptotics of the heat trace $Z(t) = \sum_{n=1}^\infty e^{-t \lambda_n}$ by Ba\~nuelos et al. \cite{Banuelos2014} and Park and Song \cite{Park2014}.
\end{abstract}

\thanks{This is a pre-print of an article published in Annales Henri Poincar\'e. The final authenticated version is available online at: https://doi.org/10.1007/s00023-018-0708-0}

\maketitle







\section{Introduction and main results}

Let $d\in\mathbb N$ and let $\Omega\subset\mathbb R^d$ be open. For $m\geq 0$, let $A_m^{\Omega}=(\sqrt{-\Delta{+}m^2}{-}\h m)_D$ denote the self-adjoint operator in $L^2(\mathbb R^d)$ defined by the closed quadratic form 
\begin{equation} \label{quad:def}
q_{m}^{\Omega}(u) = \int_{\mathbb R^d} \big(\sqrt{|2\pi\xi|^2 + m^2} - m\big)  \, |\hat u(\xi) |^2 \, d\xi
\end{equation}
with form domain \smash{$\mathcal D\big(q_m^{\Omega}\big) = H_0^{1/2}(\Omega)$}. Here, $\hat u$ denotes the Fourier transform of $u$ (with the convention of the factor $2\pi$ in the exponent), and \smash{$H_0^{1/2}(\Omega)$} denotes the fractional Sobolev space $H_0^s(\Omega)$ of order $s=\frac{1}{2}$.

If $\Omega$ is bounded, then $A_{m}^\Omega$ has compact resolvent, since the embedding \smash{$H_0^{1/2}(\Omega) \hookrightarrow$} $L^2(\Omega)$ is compact, and its spectrum consists of eigenvalues $\{\lambda_n\}_{n\in\mathbb N}$, with $0 < \lambda_1 < \lambda_2 \leq \lambda_3 \leq \dots$, accumulating at infinity only. 

With Weyl's Law \cite{Weyl1911} from 1911, it was first discovered that there is an explicit connection between geometric properties of the domain and the asymptotic growth of the counting function $N(\lambda) \defeq \sum_{n}(\lambda_n{-}\lambda)_-^0$ for the eigenvalues $\lambda_n$ of the Dirichlet Laplacian $(-\Delta)_D$, where $\np{s} \defeq -\min\{0,s\}$ for all $s\in\mathbb R$. More precisely, Weyl proved for domains with piecewise smooth boundary that the leading term in the large-$\lambda$ asymptotics of $N(\lambda)$ is proportional to the volume of the domain. Two years later, he also conjectured that the subleading term should be proportional to the surface area of the boundary. Since then, besides of trying to prove his conjectured two-term expansion, many authors have extended and improved upon Weyl's original result. While extensions to less regular domains often make use of small-time expansions of the heat trace $Z(t) = \sum_{n=1}^\infty e^{-t \lambda_n}$ (see e.g. \cite{Brown1993}) by means of the Hardy-Littlewood Tauberian theorem \cite{Hardy1914}, Ivri{\u\i}'s \cite{Ivrii1980} and Melrose's \cite{Melrose1980} two-term expansions from 1980 are based on methods from microlocal analysis and Riemannian geometry, respectively, and therefore only apply to highly regular domains. 

The fact that the asymptotic behaviour of certain functions of the eigenvalues is connected to the geometry of the domain is not a unique feature of the Dirichlet Laplacian or, more generally, of elliptic differential operators \cite{Garding1953}, but is also observed for non-local operators: In \cite{Blumenthal1960}, based on their asymptotic results on Markov operators \cite{Blumenthal1959}, Blumenthal and Getoor extend Weyl's Law to the Dirichlet fractional Laplacian $((-\Delta)^{\alpha/2})_D$, $\alpha \in (0,2]$, on Lipschitz domains. Similarly, small-time asymptotics of the heat trace of $(-\Delta)_D$ have been extended to $((-\Delta)^{\alpha/2})_D$ by Ba\~nuelos et al. \cite{Banuelos2009}, and recently to $((-\Delta{+}m^{2/\alpha})^{\alpha/2}{-}m)_D$ by Park and Song \cite{Park2014} and Ba\~nuelos et al. \cite{Banuelos2014}. More precisely, for Lipschitz domains, Park and Song prove that
\begin{equation}\label{asymp:nonlocalheattrace}
Z(t) = D^{(1)}_{\alpha} \h |\Omega| \, t^{-d/\alpha} - (D^{(2)}_{\alpha}\h|\partial\Omega| \, {-} \, m\h D^{(3)}_{\alpha}\h |\Omega|)  \,t^{-(d-1)/\alpha}+ o(t^{-(d-1)/\alpha})
\end{equation}
as $t\to 0^+$, where $D_\alpha^{(1)}, D_\alpha^{(2)}$ and $D_\alpha^{(3)}$ are positive constants only depending on $\alpha \in (0,2]$ and $d\geq 2$. For domains with $C^{1,1}$ boundary, i.e when the boundary is locally given by a differentiable function whose derivative is Lipschitz continuous, they improve the remainder to $\mathcal O(t^{-(d-2)/\alpha})$.

Asymptotic formulas for the heat trace are usually more detailed and known for more general domains than those for the counting function. As a step in between, the Riesz mean $R(\lambda)\defeq\sum_{n=1}^\infty \npb{\lambda_n{-}\lambda}$ can be obtained by integrating $N(\lambda)$, while on the other hand, $Z(t)$ can be obtained from the Laplace transform of $R(\lambda)$ (see \eqref{myheattraceresult} below). Recently, Frank and Geisinger \cite{Frank2013} obtained for the Riesz mean of the eigenvalues of $((-\Delta)^{\alpha/2})_D$ on $C^{1}$ domains the two-term asymptotic expansion
\begin{equation}
\sum_{n=1}^\infty \npb{\lambda_n{-}\lambda} \ = \ L_\alpha^{(1)} |\Omega| \, \lambda^{1+d/\alpha} - L_\alpha^{(2)} |\partial \Omega| \lambda^{1+(d-1)/\alpha} + o\big(\lambda^{1+(d-1)/\alpha}\big) \label{frank}
\end{equation}
as $\lambda\to \infty$, where $L^{(1)}_\alpha, L^{(2)}_\alpha >0$ only depend on $\alpha \in (0,2]$ and $d\geq 2$. In \cite{Frank2013}, the asymptotic formula \eqref{frank} is proved for $C^{1,\gamma}$ domains (domains with boundaries that are locally given by a differentiable function with H{\"o}lder continuous derivative with exponent $\gamma$) for any $0<\gamma\leq 1$, with a remainder whose order depends on $\gamma$. But, as noted in \cite{Frank2016}, the stated result follows for $C^1$ domains by the same argument as in \cite{Frank2012} (see also the proof of our Theorem \ref{thm:main} in Section \ref{section:proofofthm1}).

In this work, we extend the case $\alpha=1$ of \eqref{frank}, i.e. the large-$\lambda$ asymptotics of $R(\lambda)$ for the eigenvalues of \smash{$(\sqrt{-\Delta})_D$}, to $(\sqrt{-\Delta{+}m^2}{-}m)_D$ for $m>0$. The most notable difference to the massless case is the fact that 
\begin{equation}
\psi_m(t) \defeq \sqrt{t+m^2}-m \, \label{def:psi}
\end{equation} 
fails to be homogeneous in $t\geq 0$. Thus, even though the overall structure of the proof is similar to \cite{Frank2013}, the lack of homogeneity requires to use different techniques and to approach problems differently. One of the key tools we use to overcome these difficulties is the integral representation, 
\begin{equation}
q_m^\Omega(u) \, = \, \left(\frac{m}{2\pi}\right)^{(d+1)/2} \int_{\mathbb{R}^d} \int_{\mathbb{R}^d} |u(x){-}u(y)|^2 \, \frac{K_{(d+1)/2}(m|x{-}y|)}{|x{-}y|^{(d+1)/2}} \, dx\, dy \label{quad:rep}
\end{equation}
for all $u\in H_0^{1/2}(\Omega)$, where $K_\beta$ denotes the Modified Bessel Function of the Second Kind of order $\beta$.

Due to the inhomogeneity of $\psi_m$, the statement of Theorem \ref{thm:main} involves a new parameter $\mu>0$. In order to obtain an asymptotic expansion of $\sum_{n\in\mathbb N}(h\lambda_n{-}1)_-$ as $h\to 0$ for the eigenvalues $\lambda_n$ of $A_m^\Omega$, we apply Theorem \ref{thm:main} with $\mu=hm$ in Corollary \ref{thm:asymptotics} below. 

\begin{thm} \label{thm:main} For $h,\mu>0$, $d\geq 2$, and a bounded open subset $\Omega\subset \mathbb R^d$, let
\[
H_{\mu,h}^\Omega \defeq h\h A_{\mu/h}^{\Omega} - 1 =  \big(\sqrt{{-}h^2\Delta {+}\h \mu^2}\h{-}\h\mu \h{-}\h 1\big)_D 
\]
with Dirichlet boundary condition on $\Omega$. If the boundary $\partial \Omega$ belongs to $C^1$, then for all $h,\mu>0$,
\begin{equation}\label{theorem1}
\tr \big(H_{\mu,h}^\Omega\h\big)_- \, = \, \Lambda^{(1)}_{\mu} \, |\Omega| \, h^{-d} - \Lambda^{(2)}_{\mu} \, |\partial\Omega|\, h^{-d+1} + R_\mu(h) \, , 
\end{equation}
with $(1{+}\mu)^{-d/2} R_\mu(h) \in o(h^{-d+1})$ uniformly in $\mu>0$, as $h\to 0^+$. When $\partial \Omega$ belongs to $C^{1,\gamma}$ for some $\gamma>0$, then $\forall \varepsilon \in (0,\gamma/(\gamma{+}2))$, there exists $C_\varepsilon(\Omega)>0$ such that for all $h,\mu>0$, 
\[
|R_\mu(h)|\leq C_\varepsilon(\Omega) \h (1{+}\mu)^{d/2}\h h^{-d+1+\varepsilon} \, .
\] 
Here, $|\Omega|$ denotes the volume and $|\partial \Omega|$ the surface area of $\Omega$. Moreover,
\[
\Lambda^{(1)}_{\mu} \defeq \frac{1}{(2\pi)^{d}} \int_{\mathbb R^d} \Big(\sqrt{|\xi|^2{+}\mu^2} - \mu -1\Big)_- d\xi \, 
\]
and 
\[
\Lambda_\mu^{(2)}  \defeq \frac{1}{(2\pi)^{d}} \int_0^\infty \int_{\mathbb R^d} \Big(\sqrt{|\xi|^2{+}\mu^2}-\mu -1\Big)_- \Big(1{-}2F_{\mu/|\xi'|, |\xi_d|/|\xi'|}(|\xi'|t)^2 \Big) \, d\xi \, dt \, ,
\]
where $\xi'=(\xi_1,\dots,\xi_{d-1})$, and, for $\omega\geq 0$, $F_{\omega,\lambda}$ are the generalized eigenfunctions of the one-dimensional operator \smash{$((-d^2/dt^2{+}1{+}\omega^2)^{1/2}-(1{+}\omega^2)^{1/2})_D$} with Dirichlet boundary condition on the half-line, given by Kwa{\'s}nicki in \textnormal{\cite{Kwasnicki2011}} (see also Lemma \ref{lmma:Kwasnicki} in Appendix \ref{sec:halfline}).
\end{thm}

An explicit computation shows that there exists $C>0$, such that for all $\mu>0$
\begin{equation}\label{ineq:firstconstant}
\Big|\Lambda^{(1)}_{\mu} -  \Lambda^{(1)}_0 - \tfrac{\omega_d}{(2\pi)^d}  \, \mu \Big| \, \leq  \, C\, \mu^2 \, ,
\end{equation}
where \smash{$\Lambda_0^{(1)} = (2\pi)^{-d} \int_{\mathbb R^d} (|p|{-}1)_- dp =(2\pi)^{-d} (d{+}1)^{-1} \omega_d$} is the Weyl constant of $(\sqrt{-\Delta})_D$. Similarly, by a detailed analysis of the generalized eigenfunctions $F_{\omega, \lambda}$ (see Appendix F.1 in \cite{Gottwald2016}),
for $d\geq 2$ and any $\delta\in (0,1)$, there exists $C_\delta>0$, such that for all $\mu>0$, 
\begin{equation} \label{ineq:secondconstant}
\big|\Lambda^{(2)}_\mu -\Lambda^{(2)}_0 \big| \ \leq \ C_\delta \, \mu^\delta \, ,
\end{equation}
where $\Lambda^{(2)}_0 = L_1^{(2)}>0$ denotes the second constant in \eqref{frank} for $(\sqrt{-\Delta})_D$. Hence, the function \smash{$\mu\mapsto \Lambda_\mu^{(2)}$} is H\"older-continuous in $\mu=0$ and thus \smash{$\Lambda_\mu^{(2)}>0$} for $\mu$ small enough. For general $\mu$, the sign of \smash{$\Lambda_\mu^{(2)}$} remains unknown, due to the lack of a closed form expression for the generalized eigenfunctions $F_\lambda$. However, since we rely on the technical inequality \eqref{ineq:secondconstant} to obtain \eqref{thm:asymp:0} below, this does not affect the two leading orders of our asymptotic expansion.

\medskip 

By substituting $\mu=hm$ in Theorem \ref{thm:main}, we obtain from \eqref{ineq:firstconstant} and \eqref{ineq:secondconstant}:

\begin{cor} \label{thm:asymptotics} For $m\geq 0$, $d\geq 2$, $n\in\mathbb N$, and a bounded open subset $\Omega\subset \mathbb R^d$ let $\lambda_n$ denote the $n$-th eigenvalue of $A_m^\Omega = (\sqrt{-\Delta + m^2} - m)_D$. In the case when $\partial \Omega$ belongs to $C^1$, then for all $h>0$,
\begin{equation} \label{thm:asymp:0}
\sum_{n\in\mathbb N} \big(h\lambda_n-1\big)_- \ = \ \Lambda_0^{(1)} \, |\Omega| \, h^{-d} - \Big(\Lambda_0^{(2)}|\partial\Omega|{-}\,C_d \, |\Omega| \, m \Big) \, h^{-d+1} + r_m(h) \, ,
\end{equation}   
with $(m^2{+}(1{+}m)^{d/2})^{-1} r_m(h) \in o(h^{-d+1})$, uniformly in $m>0$, as $h\to 0^+$. 
When $\partial \Omega$ belongs to $C^{1,\gamma}$ for some $\gamma>0$, then $\forall \varepsilon \in (0,\gamma/(\gamma{+}2))$ there exists $C_\varepsilon(\Omega)>0$ such that 
\[
|r_m(h)|\leq C_\varepsilon(\Omega) \, \big(m^2{+}(1{+}mh)^{d/2} \big) \, h^{-d+1+\varepsilon} \, .
\]
Here, \smash{$C_d:=\frac{\omega_d}{(2\pi)^d}$}, \smash{$\Lambda_0^{(1)}=\frac{C_d}{d{+}1}$}, and \smash{$\Lambda_0^{(2)}$} is given in Theorem \ref{thm:main} and coincides with $L_1^{(2)}$, the second constant in \eqref{frank} for $\alpha = 1$.
\end{cor}

Interestingly, the sign of the subleading term in \eqref{thm:asymp:0} depends on the relation of surface area and volume of the domain $\Omega$, as well as on the value of $m$. In particular, if \smash{$|\partial \Omega|/|\Omega| > C_d \, m /\Lambda_0^{(2)}$}, the two-term asymptotics is a negative correction to the one-term Weyl asymptotics, while for \smash{$|\partial \Omega|/|\Omega| < C_d \, m /\Lambda_0^{(2)}$}, the Weyl asymptotics is exceeded. Hence, in order to keep the subleading contribution negative, for large mass $m$, the surface area needs to be much larger than the volume (and vice versa).

\medskip 
The following sections cover the proof of Theorem 1 in logical order: First, we use a continuous localization technique to be able to study the problem separately in \textit{the bulk} of the domain $\Omega$, where the boundary is not seen, and close to the boundary, where the spacial symmetries allow the reduction to a problem on the half-line.

\bigskip

\section{Localization} \label{sec:localization}
Following \cite{Solovej2003}, \cite{Solovej2010}, and \cite{Frank2013}, by constructing a continuous family of localization functions $\{\phi_u\}_{u\in\mathbb R^d}\subset C_0^1(\mathbb R^d)$ satisfying a continuous IMS-type formula (\eqref{loc:eq:IMS} below), the analysis of $\tr(H_{\mu,h}^\Omega)_-$ can be reduced to the analysis of \smash{$\tr(\phi H_{\mu,h}^\Omega \phi)_-$} with $\phi$ having support either completely contained in $\Omega$ or intersecting the boundary $\partial \Omega$. 

\begin{prop}[\textbf{Localization error}] \label{prop:loc}
For $0<l_0<\tfrac{1}{2}$, let $l:\mathbb R^d \to [0,\infty)$ be given by
\begin{equation} \label{loc:def:lu} l(u) \defeq \frac{1}{2} \Big(1+\big(\delta(u)^2 + l_0^2\big)^{-1/2}\Big)^{-1}, \end{equation}
with $\delta(u):=\mathrm{dist}(u,\Omega^c)$, where $\mathrm{dist}(u,\Omega^c) = \inf\{|x-u|: x\in\Omega^c\}$. There exists $C>0$ and $\{\phi_u\}_{u\in\mathbb R^d}\subset C_0^1(\mathbb R^d)$ s.th. $\supp \phi_u \subset B_{l(u)}(u)$, $\|\phi_u\|_\infty\leq C$, $\|\nabla \phi_u\|_\infty \leq C\, l(u)^{-1}$ for a.e. $u\in \mathbb R^d$, and
\begin{equation} \label{loc:int_phiu}  
\int_{\mathbb R^d} \phi_u(x)^2 \, l(u)^{-d} du = 1\h \quad \forall x\in\mathbb R^d \, .
\end{equation}
Moreover, there exists $C'>0$ s.th. for all $\mu>0$, $0<l_0 < \frac{1}{2}$, and $0<h\leq \frac{l_0}{8}$,
\begin{align} 0 \ & \leq \ \tr \h\npb{H_{\mu,h}^\Omega} - \int_{\mathbb R^d} \tr\h \npb{\phi_u H_{\mu,h}^\Omega \phi_u} \, l(u)^{-d}\, du \nonumber \\[3pt]
& \leq \  C'\, h^{-d+2}\h l_0^{-1} \h {\mathfrak S}_d(l_0/h) \, (1{+}\mu)^{d/2} \, , \label{ineq:prop:loc} \end{align}
where 
\begin{equation}
\mathfrak S_d(t) \, \defeq \, \left\{\begin{array}{ll}1\, , & d>2\\
|\ln(t)|^{1/2} \, , & d=2 \, .\end{array}\right. \label{defmathfrakS}
\end{equation}

\end{prop}
\begin{proof} For the construction of $\{\phi_u\}_{u\in\mathbb R^d}$ with the given properties, see \cite[Theorem 22]{Solovej2010}. In complete analogy to the case $d=3$ covered in \cite[Theorems 13, 14]{Solovej2010}, for each $u\in\mathbb R^d$, there exists a bounded linear operator $L_u$ in $L^2(\mathbb R^d)$ such that for all $f\in H_0^{1/2}(\Omega)$, 
\begin{equation} \label{loc:eq:IMS} q_{\mu/h}^\Omega(f) = \int_{\Omega^\ast} q_{\mu/h}^\Omega(f\phi_u) \, l(u)^{-d} \, du - \int_{\Omega^\ast} (f,L_{u} f) \, l(u)^{-d}\, du \, ,
\end{equation}
where $\Omega^\ast \defeq\{u\in\mathbb R^d \,|\, \supp \phi_u \cap \Omega \not= \emptyset\}$, and $\|L_u\|\leq C\,\mu^{-1} h\, l(u)^{-2}$ for all $u\in\mathbb R^d$. Moreover, analogous to \cite[Theorem 10]{Lieb1988}, it can be shown that there exists $C>0$ such that for all $0<\delta\leq\tfrac{1}{2}$ and all positive definite trace class operators $\rho$,
\begin{equation}
\tr \, \rho L_{u} \ \leq \  C \,l(u)^{-1} \Big( \delta \ \tr \, (\rho\h \chi_\Omega\h \chi_{u,\delta}) +  \tau_d(\delta)\, \|\rho\|\Big)  \qquad \forall u\in\mathbb R^d \, ,\label{loc:error} 
\end{equation}
where $\chi_{u,\delta}$ denotes the characteristic function of the ball $B_{l(u)(1+\delta)}(u)$, and 
\begin{equation}\label{taud}\nonumber
\tau_d(\delta) \, \defeq \, \left\{\begin{array}{ll}\delta^{-d+2}\, , & d>2\\
|\ln(\delta)| \, , & d=2 \, .\end{array}\right.{}
\end{equation}
For a detailed proof of these results, see \cite[Section 1]{Gottwald2016}. Similar to \cite[Proposition 1.3]{Frank2013}, from \eqref{loc:eq:IMS} and \eqref{loc:error}, it follows that
\begin{align} \nonumber
\tr \,\rho A_{\mu/h}^\Omega \ \geq & \ \int_{\Omega^\ast} \tr\,\Big(\rho \phi_u \big(A_{\mu/h}^\Omega - C \h l(u)^{-1} \h \delta_u\big )\phi_u \Big)\frac{du}{l(u)^d}  \\
&  \ - C\, \|\rho \|\int_{\Omega^\ast} \tau_d(\delta_u) \h \frac{du}{l(u)^{d+1}} \, ,  \label{rhoAgeq}
\end{align}
where
\begin{equation}\label{choiceofdeltau} \nonumber
\delta_u = \left\{ \begin{array}{lcl} l(u)^{-1} h & , & d>2 \\  l(u)^{-1} h \, |\ln(l(u)/h)|^{1/2} & , &  d=2 \, . \end{array}\right.
\end{equation}
By the Variational Principle for the sum of negative eigenvalues (see e.g. \cite{Lieb1983}), it follows from \eqref{rhoAgeq} for $d>2$,
\begin{align*} 
\tr\h \npb{H_{\mu,h}^\Omega}
& = \, - \inf_{0\leq \rho \leq 1} \tr \h \rho H_{\mu,h}^\Omega \, = \,  - \inf_{0\leq \rho \leq 1} \tr \h \rho \big(h A_{\mu/h}^\Omega - 1\big)  \\
& \leq \, \int_{\Omega^\ast} \tr\,\npb{\phi_u \big(H_{\mu,h}^\Omega- C \h h^2 l(u)^{-2}\big )\phi_u} \ \frac{du}{l(u)^d} \\
& \quad + C \h h^{-d+2} \int_{\Omega^\ast} l(u)^{-2} \, du\, .
\end{align*}
For any family $\{\sigma_u\}_{u\in\mathbb R^d}$ with $0<\sigma_u\leq \frac{1}{2}$ for all $u\in\mathbb R^d$, we have
\begin{align} \label{loc:ineq22}
&\inf_{0\leq \rho \leq \mathbb I} \tr\, \rho \big( \phi_u \big(H_{\mu,h}^\Omega - C \h h^2 l(u)^{-2}\big)\phi_u\big)\\[3pt]
& \geq (1{-}\sigma_u) \inf_{0\leq \rho \leq \mathbb I} \tr\, \rho \phi_u H_{\mu,h}^\Omega \phi_u + \inf_{0\leq \rho \leq \mathbb I} \tr\,\rho \big( \phi_u \big(\sigma_u H_{\mu,h}^\Omega - C \h h^2 l(u)^{-2}\big)\phi_u\big) \, ,  \nonumber
\end{align}
in particular, by the Variational Principle, for all $u\in\mathbb R^d$, 
\begin{align*}
&\tr\,\npb{\phi_u \big(H_{\mu,h}^\Omega- C \h h^2 l(u)^{-2}\big )\phi_u}\\[3pt]
& \leq \, \tr \, \npb{\phi_u H_{\mu,h}^\Omega \phi_u} + \tr \npb{\phi_u \big(\sigma_u H_{\mu,h}^\Omega- C \h h^2 l(u)^{-2}\big )\phi_u} \, .
\end{align*}
With $\sigma_u = h^2 l(u)^{-2}$ (note that $h^2l(u)^{-2}< \frac{1}{2}$ if $h\leq l_0/8$, since $l(u)> \frac{l_0}{3}.$), it follows that
\begin{align} \nonumber
& \tr\, \np{H_{\mu,h}^\Omega} - \int_{\Omega^\ast} \tr \,\npb{\phi_u H_{\mu,h}^\Omega \phi_u} \, \frac{du}{l(u)^d} \\ \nonumber
& \qquad \leq h^2\int_{\Omega^\ast} \tr\, \npbb{\phi_u\big(h A_{\mu/h}^\Omega - C\big)\phi_u} \ l(u)^{-d-2} du +  C \h h^{-d+2} \int_{\Omega^\ast} l(u)^{-2} \, du
\\
& \qquad \leq C \h h^{-d+2} (1{+} \mu)^{d/2}  \h  \int_{\Omega^\ast} l(u)^{-2} du\, , \label{eq04:prop:loc}
\end{align}
where we have used the inequality \eqref{bulk:lemma1:ineq} below, which yields 
\[
\tr\, \npbb{\phi\big(h A_{\mu/h}^\Omega - C\big)\phi} \, = \, C\, \tr \, \npb{\phi H_{\mu/C,h/C}^\Omega \phi} \, \leq \, C (1{+}\mu)^{d/2} h^{-d} \|\phi\|_2^2 \,
\]
for any $\phi\in C_0^1(\mathbb R^d)$. 

Let $g(t):=(l_0^2 + t^2)^{-1}$, such that $l(u)^{-2} \leq 8(1{+}g(\delta(u)))$. Thus, we want to find an upper bound for
\[
\int_{\Omega^\ast} g(\delta(u)) \h du \, = \, \int_{\delta(\Omega^\ast)} g(t) \, \big(\lambda^d\circ \delta^{-1}\big)(dt) \, ,
\]
where $\lambda^d$ denotes the Lebesgue measure on $\mathbb R^d$, and $\lambda^d\circ \delta^{-1}$ is the image measure of $\lambda^d$ under $\delta$. By the co-area formula \cite[3.4.2]{Evans1992} applied to $\delta$, we have $(\lambda^d\circ \, \delta^{-1})(dt) =  \mathcal H^{d-1}(\delta^{-1}(\{t\}))\, dt$, where $\mathcal H^{d-1}$ denotes the $(d{-}1)$-dimensional Hausdorff measure. It can be shown (see e.g. Lemma \ref{lmma:parallelsurface}) that there exist constants $\varepsilon >0$ and $C>0$ such that $\mathcal H^{d-1}(\delta^{-1}(\{t\}))\leq C$ for all $t\leq\varepsilon$. With $\Omega_\varepsilon^\ast := \{ u\in \Omega^\ast\,|\, \delta(u)\leq \varepsilon \}$, we obtain
\begin{align} \label{intdistOmega}
\int_{\Omega^\ast} g(\delta(u)) \, du \, & \leq \, C \int_0^\varepsilon g(t) \, dt + \varepsilon^{-2} \int_{\Omega^\ast{\setminus}\Omega^\ast_\varepsilon} \, du \\ & \leq \, C \, l_0^{-1} \int_0^\infty \arctan'(s) \h ds \, \leq \, C \, l_0^{-1} \, , \nonumber
\end{align}
and thus 
\begin{equation} \label{intludu}
\int_{\Omega^\ast} l(u)^{-2} du \, \leq \, C\, l_0^{-1}\, .
\end{equation}
Hence, by \eqref{eq04:prop:loc},
\begin{equation} \tr \h\npb{H_{\mu,h}^\Omega} - \int_{\mathbb R^d} \tr\h \npb{\phi_u H_{\mu,h}^\Omega \phi_u} \, l(u)^{-d}\, du \ \leq \  C\, h^{-d+2}\h l_0^{-1} (1{+}\mu)^{d/2} \, , \label{loc:ineq:prop1}\end{equation}
establishing the second inequality in \eqref{ineq:prop:loc} for $d>2$. The case $d=2$ follows along the same lines.
\end{proof}

\bigskip

\section{The bulk} 

Similar to \cite[Proposition 1.4]{Frank2013}, in the case when the support of $\phi$ is completely contained in $\Omega$, we obtain

\begin{prop}[\textbf{Error in the bulk}]\label{prop:bulk}
There exists $C>0$, such that for all real-valued $\phi\in C_0^1(\Omega)$ satisfying $\|\nabla \phi\|_\infty \leq C \h l^{-1}$ and with support in a ball of radius $l\h {>}\h 0$, 
\begin{equation}
0 \leq \Lambda_{\mu}^{(1)} \h h^{-d} \hspace{-3pt}\int_\Omega\phi(x)^2 \, dx  \, - \ \tr \npb{\phi H_{\mu,h}^\Omega\phi} \leq C \, h^{-d+2} \, l^{d-2} \, (1{+}\mu)^{(d-1)/2} \label{eq:bulk}
\end{equation}
for all $h>0$, where $\Lambda_\mu^{(1)}= (2\pi)^{-d} \int_{\mathbb R^d} (\psi_\mu(|p|^2) - 1)_- \, dp$. 
\end{prop}

\begin{proof}
We write \smash{$H_{\mu,h}:=H_{\mu,h}^{\mathbb R^d} = hA_{\mu/h}^{\mathbb R^d}{-}\h 1 = \sqrt{{-}h^2\Delta {+} \mu^2}{-}\h\mu\h {-} \h 1$} for the operator on all of $\mathbb R^d$. For the lower bound, using the the Variational Principle we obtain
\begin{equation}
\tr\np{\phi H_{\mu,h}^\Omega\phi} \, \leq \, \tr\np{\phi H_{\mu,h} \phi} \ \leq \ \tr\,\phi\np{H_{\mu,h}}\phi \, .
\label{bulk:lemma1:ineq1}
\end{equation}

If $\mathcal F$ denotes the Fourier transform on $\mathbb R^d$ and $\Phi: = \mathcal F\phi\mathcal F^{-1}$, which is the integral operator in $L^2(\mathbb R^d)$ with kernel $(k,\tilde k)\mapsto (\mathcal F\phi)(k{-}\tilde k)$, then we have $\tr\, \phi \np{H_{\mu,h}} \phi = \tr \, \Phi g \Phi$, where $g(k):= \np{\psi_\mu(|2\pi h k|^2) {-} 1}$. Since, for any $\delta\geq 0$, we have
\[ 
\int_{\mathbb R^d}\int_{\mathbb R^d} g(k)^\delta |\Phi(k,\tilde k)|^2 d\tilde k\h dk \, = \, \int_{\mathbb R^d} g(k)^\delta dk \, \|\mathcal F \phi\|_2^2 \,< \, \infty \, ,
\]
it follows that the operators $\Phi g^0$ and $g\Phi$ are Hilbert-Schmidt operators, and therefore
\begin{align*}
\tr \, \Phi g \Phi \, & = \int_{\mathbb R^d} g(k) |\Phi(k,\tilde k)|^2 \, d\tilde k \, dk \, \\ & = \int_{\mathbb R^d}\npb{\psi_\mu(|2\pi hk|^2)-1} dk \ \|\mathcal F\phi\|_2^2  \, = \, \Lambda_\mu^{(1)} \h h^{-d} \h \|\phi\|_2^2 \, .
\end{align*}
Together with \eqref{bulk:lemma1:ineq1}, this proves 
\begin{equation} 
\tr \big(\phi H_{\mu,h}^\Omega \phi\big)_- \ \leq \ \, \Lambda_\mu^{(1)} \h h^{-d} \int_{\mathbb R^d} \phi(x)^2\h dx \,   \label{bulk:lemma1:ineq}
\end{equation}
for any real-valued $\phi\in C_0^1(\mathbb R^d)$ and $h>0$.

For the upper bound, let $\rho \defeq \phi^{0} \, (H_{\mu,h})_-^0 \, \phi^0$, where $\phi^0$ is the characteristic function of $\supp \phi$. Then $0\h {\leq}\h \rho \h {\leq}\h 1$, $\rho$ is trace class, and $\phi \rho f \in H^{1/2}(\mathbb R^d)$, in particular $\rho$ can be used as a trial density matrix in the Variational Principle. It follows that
\begin{align}\label{bulk:lemma1:ineq40}
 & -\tr \, \npb{\phi H_{\mu,h}^\Omega\phi} \,  \leq \, \tr \, \rho \phi H_{\mu,h}\phi \, = \, h \tr\, \rho \phi A_{\mu/h} \phi -  \tr \,\rho \phi^2 . \\ \nonumber
& = \, (2\pi h)^{-d} \int_{\psi_\mu(|p|^2)\leq 1}  \Big(\big\|(hA_{\mu/h})^{1/2} \phi\h e^{ip\, \cdot /h} \big \|^2_2 -  \big\|\phi e^{ip\cdot/h}\big\|_2^2 \Big) \,  dp  \, . 
\end{align}
By using an exponential regularization for $\psi_\mu$, $\psi_\mu^a(t):= e^{-aE} \psi_\mu(t)$ for all $a>0$ and $t>0$, and writing $\phi(x)\phi(y) = \tfrac{1}{2} \big(\phi(x)^2 + \phi(y)^2 - (\phi(x) {-} \phi(y))^2\big)$, the first term in \eqref{bulk:lemma1:ineq40} can be written as
\begin{align}\nonumber
&\big\|(h\h A_{\mu/h})^{{1}/{2}} \phi \h e^{ip\cdot/h}\big\|_2^2\\
& \  = \, \frac{1}{2} \int_{\mathbb R^d} \big(\psi_\mu(|p{+}2\pi h\eta|^2) + \psi_\mu(|p{-}2\pi h\eta|^2)\big)\, |\hat \phi(\eta)|^2 \, d\eta \, . \label{bulk:eq50}
\end{align}
We omit the proof of \eqref{bulk:eq50}, since it is purely technical (see \cite[Lemma 10]{Gottwald2016}). From \eqref{bulk:lemma1:ineq40} and \eqref{bulk:eq50}, we obtain
\begin{align}
-\tr \, \npb{\phi H_{\mu,h}^\Omega \phi} & 
\, \leq \, (2\pi h)^{-d} \int_{\psi_\mu(|p|^2)\leq 1} \Big(\big\|(hA_{\mu/h})^{1/2} \phi\h e^{ip\, \cdot /h} \big \|^2_2 - \|\phi\|_2^2\Big) \,  dp \, \nonumber \\
& \, = \, - h^{-d} \Lambda_\mu^{(1)} \|\phi\|_2^2 + (2\pi h)^{-d} \int_{\psi_\mu(|p|^2)\leq 1} R_{\mu,h}(p)\, dp \, , \label{eq:bulkmain4}  
\end{align}
where
\[ 
R_{\mu,h}(p) \ \defeq \ \int_{\mathbb R^d} \Big(\tfrac{1}{2} \big(\psi_\mu(|p{+}2\pi h\eta|^2) + \psi_\mu(|p{-}2\pi h\eta|^2)\big) - \psi_\mu(|p|^2) \Big) \, |\hat\phi(\eta)|^2 \, d\eta \, .
\] 
It remains to find a suitable upper bound for $R_{\mu,h}(p)$. For $a>0$ and $|b|\leq a$, we have $(a+b)^{1/2} + (a-b)^{1/2} \leq 2\h a^{1/2}$. Applied to $a= |p|^2 +|\xi|^2 +\mu^2$ and $b= 2p\cdot \xi $, where $p,\xi \in\mathbb R^d$, we find
\begin{align*}
&\tfrac{1}{2} \big(\psi_\mu(|p{+}\xi|^2) {-} \psi_\mu(|p{-}\xi|^2)\big) -\psi_\mu(|p|^2) \\
&\leq \,  \big(|p|^2{+}|\xi|^2{+}\mu^2 \big)^{1/2} {-}\, \big( |p|^2 {+}\mu^2 \big)^{1/2}  \, \leq \, \tfrac{1}{2}\, |p|^{-1} |\xi|^2 \, ,
\end{align*}
where the last inequality is due to \mbox{$(c+d)^{1/2}- c^{1/2}\leq \tfrac{1}{2} c^{-1/2} d$}, which holds for all $c,d>0$. Hence,
\begin{align}
&\int_{\psi_\mu(|p|^2)\leq 1} R_{\mu,h}(p) \, dp  
 \, \leq \, \tfrac{1}{2} (2\pi h)^2 \int_{\psi_\mu(|p|^2)\leq 1} |p|^{-1}  dp \int_{\mathbb R^d} |\eta|^2 |\hat\phi(\eta)|^2 d\eta \nonumber \\
& \, = \,  2\pi^2 \h h^2 \h \ |\mathbb S^{d-1}| \int_{0}^{\sqrt{1{+}2\mu}} \hspace{-5pt} t^{d-2}\, dt \  \|\nabla \phi \|_2^2   \, \leq \, C \h l^{d-2}\h h^2\h (1{+}\mu)^{(d-1)/2} \, , \nonumber
\end{align}
where we have used the assumption $\|\nabla \phi\|_\infty\leq C l^{-1}$ and $|\supp (\nabla \phi)|\leq Cl^{d}$. Together with \eqref{eq:bulkmain4} this proves the upper bound in \eqref{eq:bulk}.
\end{proof}

\bigskip

\section{Straightening of the boundary} \label{sec:straightening}

In this section, we compare $H_{\mu,h}^\Omega$ locally near the boundary with \smash{$H_{\mu,h}^+$}, where
\[
H_{\mu,h}^+ \defeq H_{\mu,h}^{\mathbb R^d_+} \, , \quad  \mathbb R^d_+ \ \defeq \ \big\{(\xi',\xi_d)\in \mathbb R^{d-1}{\times}\mathbb R\, \big|\,  \xi_d>0 \big\} \, .
\]
Proposition \ref{prop:straightening} below is an extension of \cite[Lemma 4.2]{Frank2013} and relies on the assumption that the boundary is locally given by the graph of a differentiable function. 

More precisely, if the support of $\phi\in C_0^1(\mathbb R^d)$ is contained in an open ball $B_l\subset \mathbb R^d$ of radius $0<l\leq c$ with $B_l\cap \partial \Omega \not = 0$ and some $c>0$ to be fixed later, we choose new coordinates in $\mathbb R^d$ in the following way: By translation and rotation, we can choose a Cartesian coordinate sytem centered at some $x_l\in B_l\cap \partial \Omega$, such that $(0,1)$ is the unit inward normal vector at $x_l=(0,0)$, where for $x\in\mathbb R^d$, we write $x=(x',x_d)\in\mathbb R^{d-1} \h{\times}\h \mathbb R$. Let
\[
D_l  \ \defeq \ \big\{x'\in\mathbb R^{d-1} : (x',x_d)\in B_l \big\}
\] 
be the projection of $B_l$ on the hyperplane $\partial\mathbb R^d_+ \defeq \{(x',x_d)\in \mathbb R^d : x_d \,{=}\, 0\}$. Since $\partial \Omega \in C^{1}$, for small enough $c>0$ and $l\leq c$, there exists a differentiable function $g:D_l \to \mathbb R$, such that
\[ 
B_l \cap \partial \Omega = \big\{(x',g(x')): x'\in D_l \big\}\, .
\]
Moreover, since $\partial \Omega$ is compact, the derivatives of the functions $g$ for different patches along $\partial \Omega$ admit a common modulus of continuity $w:\mathbb R_+\to\mathbb R_+$. In particular, $w$ is non-decreasing, $w(t)\to 0$ as $t\to 0^+$, and 
\begin{equation} \label{modulusofcontinuity}
|\nabla g(x'){-}\nabla g(y')|\leq w(|x'{-}y'|)
\end{equation}
for all $x',y'\in D_l$. We define the diffeomorphism
\begin{equation} 
\tau: D_l \times \mathbb R \to D_l \times \mathbb R, ~(x',x_d) \mapsto \big(x',x_d\h{-}\h g(x')\big). \label{def:T}
\end{equation}
Then $\tau$ straightens the part of $\partial \Omega$ that lies inside of $B_l$, in the sense that it maps $B_l\cap \partial \Omega$ into $\partial \mathbb R_+^d$, since $\tau(x',g(x'))=(x',0)$ for all $x'\in D_l$. 

\vspace{5pt}

\begin{prop}[\textbf{Straightening of the boundary}]\label{prop:straightening} There exist positive constants $c$ and $C$, such that for any \mbox{$\phi\in C_0^1(\mathbb R^d)$} with support in a ball of radius $0<l\leq c$ intersecting the boundary of $\Omega$, we have
\begin{equation} 
\big| \tr \, \npb{\phi H_{\mu,h}^\Omega\phi} - \tr \, \npb{\phi' H_{\mu,h}^+ \phi'} \big| \leq C \h w(l)\, l^d \h (1{+}\mu)^{d/2} h^{-d} \, . \label{eq:prop:straightening}
\end{equation}
Here, $\phi'\in C_0^1(\mathbb R^d)$ denotes the extension of $\phi\circ \tau^{-1}$ by zero to $\mathbb R^d$, where $\tau$ is the diffeomorphism given in \eqref{def:T}. Moreover,
\begin{equation}\int_{\Omega}\phi(x)^2\, dx \, = \, \int_{\mathbb R^d_+} \phi'(x)^2 \, dx \, , \label{prop:straighten:eq02}\end{equation}
and there exists $C>0$ such that
\begin{equation} 0 \, \leq \, \int_{\partial \Omega} \phi(x)^2 d\sigma(x) - \int_{\mathbb R^{d-1}} \phi'(x',0)^2 dx' \, \leq \, C \, w(l)^2 \, l^{d-1} \, . \label{prop:straighten:eq03}\end{equation}
\end{prop}

For $\mu=0$, \eqref{eq:prop:straightening} is proved in \cite[Lemma 4.2]{Frank2013} by directly using the homogeneity of $|\cdot|$. Instead, since $\psi_\mu$ is not homogeneous for $\mu>0$, we are using the integral representations \eqref{bessel:intrep2} and \eqref{bessel:intrep1} (see Appendix \ref{app:bessel}) of the modified Bessel functions in the expression \eqref{quad:rep} for the quadratic form $q_{\nu}^\Omega$ defined in \eqref{quad:def}.

\begin{lem}\label{lmma:straight1}
There exist $c, C>0$, such that for all $v\in H_0^{1/2}(\Omega)$ that are compactly supported in a ball $B_l$ of radius $0<l\leq c$, and for all $\nu>0$,
\begin{equation}
\big|q_\nu^\Omega(v)-q_\nu^+(v')\big| \leq C \h w(l) \h \min\big\{q_\nu^\Omega(v),q_\nu^+(v')\big\}\, ,\label{eq:lmma:straight1}
\end{equation}
where $q_\nu^+:=q_\nu^{\mathbb R^d_+}$, and $v'\in H_0^{1/2}(\mathbb R^d_+)$ denotes the extension of the function $v\circ \tau^{-1}$ by zero to $\mathbb R^d_+$.
\end{lem}

\begin{proof}
We start with an upper bound for the left side of \eqref{eq:lmma:straight1} in terms of $q_\nu^+(v')$. Writing $\theta_\nu(t) \defeq \nu^{d+1}\theta(\nu t) = (\nu/(2\pi t))^{(d+1)/2} \, K_{(d{+}1)/2}(\nu t)$, then, by the integral representation \eqref{quad:rep}, 
\begin{align}\nonumber
q_\nu^\Omega(v) \,
&= \, \int_{\mathbb R^d}\int_{\mathbb R^d} \big|v(x)\h{-}\h v(y)\big|^2 \, \theta_\nu(|x\h{-}\h y|) \, dx \h dy \\ \nonumber
&= \, \int_{\Gamma_l}\int_{\Gamma_l} \big|v'(\xi)\h{-}\h v'(\eta)\big|^2 \, \theta_\nu(|\tau^{-1}(\xi)\h{-}\h \tau^{-1}(\eta)|) \, d\xi \h d\eta \\
& \quad + 2 \int_{\mathbb R^d\setminus \Gamma_l}\int_{\Gamma_l}|v'(\xi)|^2 \, \theta_\nu(|\tau^{-1}(\xi)\h{-}\h y|)  \, d\xi \h dy \, , \label{eq:straight1_1}
\end{align}
since $\supp v \subset \Gamma_l\defeq D_l\times\mathbb R$, and $\tau$ has unit Jacobian determinant and is bijective on $\Gamma_l$. From the integral representation \eqref{bessel:intrep2} of $K_{(d{+}1)/2}(\nu t)$, it follows for all $\xi,\eta \in \Gamma_l$ that
\begin{align} \label{eq:straight1_2}
& \Big| \theta_\nu\big(|\tau^{-1}(\xi){-}\tau^{-1}(\eta)|\big) - \theta_\nu(|\xi{-}\eta|) \Big| \\ 
& \leq \nu^{-(d+1)/2} \, \int_0^\infty  \frac{e^{-\nu^2 u}}{(2u)^{(d+3)/2}} \Big |e^{-|\tau^{-1}(\xi)-\tau^{-1}(\eta)|^2 /(4u)} - e^{-|\xi-\eta|^2 /(4u)} \Big|  \, du \, . \quad \nonumber 
\end{align}
Since $\tau^{-1}(\xi) = (\xi',\xi_d + g(\xi'))$ for all $\xi \in \Gamma_l$ we have 
\begin{align}
\big|\tau^{-1}(\xi) {-} \tau^{-1}(\eta)\big|^2 - |\xi{-}\eta|^2 
&\nonumber  \, = \,2(\xi_d{-}\eta_d)\big(g(\xi'){-}g(\eta')\big)+ |g(\xi'){-}g(\eta')|^2 \\
& \leq 2|\xi_d{-}\eta_d|\h|\xi'{-}\eta'| \h \|\nabla g\|_\infty + |\xi'{-}\eta'|^2 \h \|\nabla g\|_\infty^2 \nonumber\\
& \leq C\, w(l) \, |\xi{-}\eta|^2 \, , \label{eq:straight1_3}
\end{align}
where we have used that, by \eqref{modulusofcontinuity}, 
\begin{equation} 
|\nabla g(x')| \, = \,|\nabla g(x'){-}\nabla g(0)| \leq  w(|x'|) \leq  w(l) \, \label{lmma:straighten:ineq0003}
\end{equation}
for all $x'\in D_l$, and therefore $\|\nabla g\|_\infty \leq w(l) <1$ for $l$ small enough.

Since $1-e^{-|t|} \leq |t|$ for all $t\in\mathbb R$, we obtain from \eqref{eq:straight1_3} that \eqref{eq:straight1_2} is bounded by
\[
C\h  w(l) \h |\xi{-}\eta|^2 \h \nu^{-(d+1)/2} \int_0^\infty  \frac{e^{-\nu^2 u-|\xi-\eta|^2/(4u)}}{(2u)^{(d+5)/2}} \, du  \, .
\]
By using the integral representation \eqref{bessel:intrep2} again, we conclude that
\begin{align}
\Big| \theta_\nu\big(|\tau^{-1}(\xi){-}\tau^{-1}(\eta)|\big) - \theta_\nu(|\xi{-}\eta|) \Big|
& \, \leq \, C \h w(l) \h \nu \h |\xi{-}\eta| \, \frac{K_{(d+3)/2}(\nu|\xi{-}\eta|)}{|\xi{-}\eta|^{(d+1)/2}} \nonumber \\
& \leq C \h w(l) \, \theta_\nu(|\xi{-}\eta|/\text{\small$\sqrt{2}$}) \, , \label{eq:straight1_4}
\end{align}
where we have used Lemma \ref{lmma:besselineq} below (based on the integral representation \eqref{bessel:intrep1} in Appendix \ref{app:bessel}), by which it follows that 
\[
\nu \h |\xi{-}\eta| \, K_{(d+3)/2}(\nu|\xi{-}\eta|)  \leq  2\h K_{(d+1)/2}\big(\nu|\xi{-}\eta| /\text{\small$\sqrt{2}$}\big) \, .
\]

Next, considering the second term in \eqref{eq:straight1_1}, containing $\theta_\nu(|\tau^{-1}(\xi){-}y|)$ with \mbox{$\xi \in \Gamma_l$} and \mbox{$y\in \mathbb R^d{\setminus}\Gamma_l$}, we have
\begin{align*}
&\Big| \theta_\nu\big(|\tau^{-1}(\xi){-}y|\big) - \theta_\nu(|\xi{-}y|) \Big| \\
& \ \leq \, \nu^{-(d+1)/2} \, \int_0^\infty  \frac{e^{-\nu^2 u}}{(2u)^{(d+3)/2}} \Big |e^{-|\tau^{-1}(\xi)-y|^2 /(4u)} - e^{-|\xi-y|^2 /(4u)} \Big|  \, du \, .
\end{align*}
As above, from $\tau^{-1}(\xi) = (\xi', \xi_d + g(\xi'))$ and $y=(y',y_d)$, it follows that
\begin{align*}
|\tau^{-1}(\xi) {-} y|^2 - |\xi{-}y|^2 \, \leq \, w(l) \, |\xi{-}y|^2\, \, ,
\end{align*} 
and therefore
\begin{equation} 
\Big| \theta_\nu \big(|\tau^{-1}(\xi){-}y|\big) - \theta_\nu\big(|\xi{-}y|\big) \Big|  \leq C \h w(l) \, \theta_\nu \big(|\xi{-}y|/\text{\small$\sqrt{2}$}\big) \, . \label{eq:straight1_6}
\end{equation}
Together with \eqref{eq:straight1_4} and \eqref{eq:straight1_1}, we obtain \smash{$| q_\nu^\Omega(v) {-} q_\nu^+(v')| \leq  C w(l) \h q_{\nu/\text{\tiny$\sqrt{2}$}}^+(v')$}. Since, for all $f\in H_0^{1/2}(\mathbb R^d_+)$,
\begin{align*}
q_{\nu/\text{\tiny$\sqrt{2}$}}^+(f) & \, = \,\tfrac{1}{\sqrt{2}} \int_{\mathbb R^d} \psi_\nu(2|2\pi k|^2) \, |\hat f(k)|^2 \, dk  \ \leq \  \sqrt{2} \, q_\nu^+(f) \, ,
\end{align*}
this finishes the proof of the first inequality in \eqref{eq:lmma:straight1}. By interchanging the roles of $q_{\nu}^\Omega$ and $q_{\nu}^+$, following the same lines as above leads to the other inequality in \eqref{eq:lmma:straight1} and thus finishes the proof of Lemma \ref{lmma:straight1}.  
\end{proof}

\vspace{10pt}
\begin{proof}[Proof of Proposition \ref{prop:straightening}] Similar to \cite[Lemma 4.2]{Frank2013}, equations \eqref{prop:straighten:eq02} and \eqref{prop:straighten:eq03} immediately follow from a change of variables and \eqref{lmma:straighten:ineq0003}. By the Variational Principle, it follows from Lemma \ref{lmma:straight1} that
\begin{align*}
\tr\h \np{\phi H_{\mu,h}^\Omega \phi}  \, \leq \, \tr \, \npbb{\phi' \Big( h\big(1{-}Cw(l) \big)A_{\mu/h}^+ -1 \Big)\phi'} \ \, . \label{ineq:prop:straight2}
\end{align*}
Moreover, as in \eqref{loc:ineq22}, for any $0<\varepsilon\leq 1/2$,
\begin{align*}
&\tr \npbb{\phi' \Big( h\big(1-Cw(l) \big)A_{\mu/h}^+ -1 \Big)\phi'}  \\
 & \qquad \qquad  \leq \, \tr \npbb{\phi'(hA_{\mu/h}^+ {-} 1)\phi'} + \tr \npbb{\phi'\big((\varepsilon {-} Cw(l)) hA_{\mu/h}^+ - \varepsilon \big)\phi'} \, .
\end{align*}
Hence, for $c$ small enough such that $\varepsilon \defeq 2 Cw(l) \leq 1/2$ for all $0<l\leq c$, it follows that 
\begin{align*}
\tr \h \npb{\phi H_{\mu,h}^\Omega \phi} - \tr \npb{\phi' H_{\mu,h}^+ \phi'}\, 
& \leq \, 2 C \h w(l) \,  \tr \h \npbb{\phi'\big(\tfrac{h}{2} A_{\mu/h}^+ {-}1\big) \phi'}  \\
& \leq \,  C \h w(l) \h l^d \h(1{+}\mu)^{d/2} h^{-d} \, ,
\end{align*}
where the second inequality is the analogue of \eqref{bulk:lemma1:ineq} for \smash{$A_{\mu/h}^+$}. Due to the symmetry of Lemma \ref{lmma:straight1} with respect to $q_\nu^\Omega$ and $q_\nu^+$, interchanging the roles of $H_{\mu,h}^\Omega$ and $H_{\mu,h}^+$ in the proof above yields the same result and therefore finishes the proof of Proposition \ref{prop:straightening}.
\end{proof}

\bigskip

\section{Analysis near the boundary}

The results of the previous section reduce the analysis of \smash{$\tr\np{\phi H_{\mu/h}^\Omega \phi}$} for $\mathrm{supp}(\phi)$ intersecting the boundary, to a problem on the half-space $\mathbb R^d_+$. In this section, it is further reduced to a problem on the half-line (Lemma \ref{hsptohl:lmma2} below). Following \cite[Section 3.2]{Frank2013}, we define a unitary operator from $L^2(\mathbb R^d_+)$ to the constant fiber direct integral space \smash{$\int^{\oplus}_{\mathbb R^{d-1}} L^2(\mathbb R_+)\defeq L^2(\mathbb R^{d-1}; L^2(\mathbb R_+))$} \cite[XIII.16]{Reed1978}. This allows to express $\smash{A_{\mu/h}^+}$ in terms of a family of one-dimensional model operators $\{\smash{T^+_{\omega}}\}_{\omega\geq 0}$, for which we apply the diagonalization results by Kwa\'snicki \cite{Kwasnicki2011} in Lemma \ref{lmma:diag} below. 

\bigskip
The main result of this section is

\begin{prop}[\textbf{Error in the half-space}] \label{prop:asymphalf}
For all $\delta_1,\delta_2\in (0,1)$ there exist constants $C_{\delta_1},C_{\delta_1,\delta_2}>0$ such that for all real-valued $\phi\in C_0^1(\mathbb R^d)$ supported in a ball of radius 1,
\begin{align} \nonumber
& - C_{\delta_1,\delta_2} \,\Big( (1{+}\mu)^{(d-\delta_2)/2} h^{-d+1+\delta_2} + (1{+}\mu)^{(d-\delta_1)/2} \, h^{-d+1+\delta_1}\Big) \\ \nonumber
& \leq \, \tr \, \big(\phi H^+_{\mu,h}\phi\big)_- - h^{-d}\,  \Lambda^{(1)}_\mu \int_{\mathbb R^d_+}\phi(x)^2 dx + h^{-d+1} \Lambda_\mu^{(2)}\int_{\mathbb R^{d-1}} \phi(x',0)^2 \, dx' \\ 
& \leq \, C_{\delta_1} \, (1{+}\mu)^{(d-\delta_1)/2} \, h^{-d+1+\delta_1} \, ,
\end{align}
where  $\Lambda_\mu^{(2)} = \int_0^\infty \mathcal K_\mu(t)\, dt$, and for $\nu,t >0$, 
\begin{align}
\mathcal K_{\mu}(t) & \defeq \, \frac{1}{(2\pi)^{d-1}} \int_{\mathbb R^{d-1}} |\xi'|^2\, \Big(\mathcal J_{\mu,|\xi'|}- \mathcal J^+_{\mu,|\xi'|}(|\xi'| t) \Big) \h d\xi' \, , \label{defKmu} \\ \nonumber
\mathcal J^+_{\mu,\nu}(t)\, & \defeq \,\frac{2}{\pi} \int_0^\infty \npb{\psi_{\mu/\nu}(\lambda^2{+}1)-\nu^{-1}} \, F_{\mu/\nu,\lambda}(t)^2 \h d\lambda \, , \\
\mathcal J_{\mu,\nu} \, & \defeq \, \frac{1}{\pi}\int_0^\infty \npb{\psi_{\mu/\nu}(\lambda^2{+}1)-\nu^{-1}} \, d\lambda \, , \nonumber
\end{align}
with $\psi_\omega$ as defined in \eqref{def:psi}.
\end{prop}

This is a generalization of \cite[Proposition 3.1]{Frank2013} for $\mu>0$. By developing an explicit diagonalization of \smash{$hA_{\mu/h}^+$}, the following two lemmas (Lemma \ref{hsptohl:lmma2} and \ref{lmma:diag}) set the basis for its proof.

\begin{lem}[\textbf{Reduction to the half-line}] \label{hsptohl:lmma2} For $\omega\h {\geq}\h 0$, let $Q_\omega^+$ denote the closed quadratic form with domain $\smash{H_0^{1/2}(\mathbb R_+)}$ and
\begin{equation} 
Q_{\omega}^+(u) \, \defeq \,  \int_{\mathbb R}  \psi_\omega\big((2\pi s)^2{+}1\big) \, |\hat u(s)|^2 \, ds \, .\label{def:qalpha}
\end{equation}
There exists a unitary operator \smash{$U:L^2(\mathbb R^d_+)\to \int^{\oplus}_{\mathbb R^{d-1}} L^2(\mathbb R_+)$} such that 
\begin{equation}
q_{\mu/h}^+(f) = \int_{\mathbb R^{d-1}} |2\pi \xi'|\, Q_{\mu/|2\pi h \xi'|}^+\Big( (Uf)_{\xi'} \Big) \, d\xi'\, \qquad \forall f\in H_0^{1/2}(\mathbb R^d_+) \, . \label{hsptohl:lmma2:eq2}
\end{equation}
\end{lem}

\begin{proof}
Let $\mathcal F^{(d-1)}$ denote the partial Fourier transform of $f\in L^2(\mathbb R^d)$ in the first $d{-}1$ variables, and for $g\in \int^{\oplus}_{\mathbb R^{d-1}} L^2(\mathbb R_+)$ and $\xi'\in \mathbb R^{d-1}$, we use the notation $g_{\xi'}\defeq g(\xi')\in L^2(\mathbb R_+)$. By the unitarity of the Fourier transform, and a change of variables, it follows that the operator $U:L^2(\mathbb R^d_+)\to \int^{\oplus}_{\mathbb R^{d-1}} L^2(\mathbb R_+)$, given by
\begin{equation}
(Uf)_{\xi'}(t) \ \defeq \ |2\pi\xi'|^{-1/2} \big( \mathcal F^{(d-1)} f\big)(\xi',|2\pi\xi'|^{-1} t)\,  \label{bdr:defU}
\end{equation}
for almost all $\xi'\in \mathbb R^{d-1}$ and $t>0$, is unitary. Since, for $\omega,\nu,s,t>0$ and $u\in L^2(\mathbb R)$, we have \smash{$\psi_\omega(\nu^2 t)= \nu \,\psi_{\omega/\nu}(t)$}, and $\big(\mathcal Fu(\omega^{-1}\cdot)\big)(s) = \omega \h (\mathcal F u)(\omega s)$, where $\mathcal F$ denotes the Fourier transform in $L^2(\mathbb R)$, equation \eqref{hsptohl:lmma2:eq2} follows from a change of variables.
\end{proof}

\bigskip 

The self-adjoint operators $T_\omega^+$ in $L^2(\mathbb R_+)$ given by the closed quadratic form $Q_\omega^+$, i.e. \smash{$T_\omega^+  = \psi_\omega \big({-}\frac{d^2}{dt^2} {+}1\big)$} with Dirichlet boundary condition on $\mathbb R_+$, belong to a large class of operators for which an explicit diagonalization in terms of generalized eigenfunctions has been proved by Kwa\'snicki in \cite{Kwasnicki2011}. In Appendix \ref{sec:halfline}, we apply these results (Corollary \ref{corspec}) and derive properties for the terms in this spectral decomposition (Lemmas \ref{hspthl:lmma:propvartheta}, \ref{hspthl:lmma:propvarphi}). 

As a consequence of Lemma \ref{hsptohl:lmma2} and Corollary \ref{corspec}, we obtain

\begin{lem}[\textbf{Diagonalization of $A^+_{\mu,h}$}] \label{lmma:diag} For $\omega\geq 0$ let $\{F_{\omega,\lambda}\}_{\lambda>0}$ denote the generalized eigenfunctions of $T_\omega^+$ given in Corollary \ref{corspec} below. The linear map $V_h$, defined on $L^1\cap L^2(\mathbb R^d_+)$ by 
\[
V_hf(\xi) = \int_{\mathbb R^d_+} v_h(\xi,x) f(x) \, dx \, \qquad \forall \xi=(\xi',\xi_d) \in \mathbb R^d_+\, ,
\]
where $v_h(\xi,x)\defeq h^{-d/2}\h v(\xi,h^{-1}x)$ and 
\begin{equation}\label{kernelofunitaryequiv}
v(\xi,x) \defeq |\xi'|^{1/2} \, \frac{e^{-i\xi'x'}}{(2\pi)^{(d{-}1)/2}}\, \sqrt{\frac{2}{\pi}} \, F_{\mu/|\xi'|,\xi_d}(|\xi'|x_d) \, ,
\end{equation}
extends to a unitary operator $V_h:L^2(\mathbb R^d) \to L^2(\mathbb R^d)$ and establishes the unitary equivalence of \smash{$hA_{\mu/h}^+$} with multiplication by \smash{$a_\mu(\xi',\xi_d):=|\xi'| \psi_{\mu/|\xi'|}(\xi_d^2{+}1)$}, i.e. 
\begin{equation}
V_hhA_{\mu/h}^+V_h^\ast = a_\mu \, . \label{unitaryequivalence}
\end{equation}
\end{lem}
\begin{proof}
For $h>0$, let $S_h$ be the unitary scaling operator in \smash{$\int^\oplus_{\mathbb R^{d-1}} L^2(\mathbb R_+)$} given by 
\[
(S_hg)_{\xi'}(t) = (2\pi h)^{-(d-1)/2} g_{\xi'/2\pi h}(t) \, ,
\] 
and let $U_h:=S_h\circ U$, with the unitary operator $U$ given in \eqref{bdr:defU}. Then, by \eqref{hsptohl:lmma2:eq2},
\begin{equation}
h\h q_{\mu/h}^+(f) \ = \ \int_{\mathbb R^{d-1}} |\xi'| \, Q_{\mu/|\xi'|}^+\Big( (U_h f)_{\xi'} \Big) \, d\xi' \, . \label{hsptohl:lmma2:eq2adv}
\end{equation}
By Corollary \ref{corspec}, for any $\omega\geq 0$, $T_\omega^+$ is unitarily equivalent to the operator of multiplication in $L^2(\mathbb R_+)$ by $\lambda\mapsto \psi_\omega(\lambda^2{+}1)$, where, for $\phi\in L^1(\mathbb R_+)\cap L^2(\mathbb R_+)$, the corresponding unitary transformation is explicitly given by
\[
\mathit{\Pi_\omega} \phi (\lambda) \ = \ \sqrt{\frac{2}{\pi}} \int_0^\infty F_{\omega,\lambda}(t) \h \phi(t)\h dt \, .
\]
It follows that
\begin{equation} \label{hsptohl:lmma2:eq3}
h\h q_{\mu/h}^+(f) \, = \, \int_{\mathbb R^{d-1}} \int_{\mathbb R_+} a_\mu(\xi',\xi_d) \, \big|\mathit{\Pi_{\mu/|\xi'|}} (U_hf)_{\xi'} (\xi_d) \big|^2 \, d\xi_d \, d\xi' \, .
\end{equation}
Moreover, since $\mathit{\Pi_{\mu/|\xi'|}} (U_hf)_{\xi'}(\xi_d) = (V_h f)(\xi',\xi_d)$ for all $\xi'\in\mathbb R^{d-1}, \xi_d\in \mathbb R_+$, after a change of variables, the unitarity of $V_h$ follows from the unitarity of the partial Fourier transform and the unitarity of $\mathit{\Pi}_{\omega}$. In particular, \eqref{unitaryequivalence} follows from \eqref{hsptohl:lmma2:eq3}.
\end{proof}

\bigskip

Using this, we derive the following representation of $\tr \, \phi (H_{\mu,h}^+)_- \phi$, which will be used to prove the lower bound in Proposition \ref{prop:asymphalf}.

\begin{lem} \label{lmma:hspasymp:01}
For any real-valued $\phi \in C_0^1(\mathbb R^d)$ we have
\begin{align} \label{asymhs:eq3}
& \tr \,  \phi \npb{H_{\mu,h}^+} \phi \ 
 = \ \int_{\mathbb R^d_+} \int_{\mathbb R^d_+}|\xi'|\big(a_\mu(\xi){-}1\big)_- \, |v_h(\xi,x)|^2 \, d\xi \, \phi(x)^2 \, dx \\
& = \ h^{-d}\,  \Lambda^{(1)}_\mu \int_{\mathbb R^d_+}\phi(x)^2 dx - h^{-d+1}\int_{\mathbb R^d_+}\phi(x)^2 h^{-1}\mathcal K_\mu(h^{-1}x_d) \, dx \, , \label{asymhs:eq4}
\end{align}
where $\mathcal K_\mu$ was defined in \eqref{defKmu}.
\end{lem}

\begin{proof} 
Since, by the definitions of $a_\mu(\xi)$, $v_h(\xi,x)$ and $\mathcal J_{\mu,|\xi'|}^+(t)$, 
\begin{align*}
\int_{\mathbb R^d_+}|\xi'|\big(a_\mu(\xi){-}1\big)_- \, |v_h(\xi,x)|^2 \, d\xi = \frac{h^{-d}}{(2\pi)^{d-1}} \int_{\mathbb R^{d-1}} |\xi'|^2\, \mathcal J^+_{\mu,|\xi'|}(h^{-1}|\xi'|x_d)\, d\xi' \, ,
\end{align*}
and, by changing variables,
\[
\frac{1}{(2\pi)^{d-1}} \int_{\mathbb R^d_+}\int_{\mathbb R^{d-1}} |\xi'|^2\, \mathcal J_{\mu,|\xi'|}\, d\xi' \, \phi(x)^2 \, dx \, = \, \Lambda_\mu^{(1)} \int_{\mathbb R^d_+} \phi(x)^2\, dx \, ,
\]
it follows from the definition of $\mathcal K_\mu$ in \eqref{defKmu}, that \eqref{asymhs:eq4} is a direct consequence of \eqref{asymhs:eq3}.

For simplicity, we write $a\defeq a_\mu$ and $V\defeq V_h$. First, we show that $(a{-}1)^0_-V\phi V^\ast$ and $(a{-}1)_- V\phi V^\ast$ are Hilbert-Schmidt operators. 

For any $0\leq \delta \leq 1$ we have
\begin{align}
& \int_{\mathbb R^d_+} \big(a(\xi){-}1\big)_-^\delta \int_{\mathbb R^d_+} |V\phi V^\ast(\xi,\zeta)|^2 d\zeta \h d\xi \nonumber \\
& = \, \int_{\mathbb R^d_+} \big(a(\xi){-}1\big)_-^\delta \int_{\mathbb R^d_+} \bigg|\int_{\mathbb R^d_+} v_h(\xi,x) \phi(x) \overline{v_h(\zeta,x)} \, dx \bigg|^2 d\zeta \h d\xi \nonumber \\
& = \lim_{c\to 0^+} \lim_{b\to 0^+} \int_{\mathbb R^d_+} \big(a(\xi){-}1\big)_-^\delta \int_{\mathbb R^d_+} e^{-c|\xi'-\zeta'|^2} e^{-bf_{\mu/|\xi'|}(\zeta_d^2)} \nonumber \\
& \ \ \qquad \qquad \qquad \qquad  \times  \bigg|\int_{\mathbb R^d_+} v_h(\xi,x) \phi(x) \overline{v_h(\zeta,x)} \, dx \bigg|^2 d\zeta \h d\xi \label{asymhs:eq1} \, ,
\end{align}
where $f_w(t)= \psi_w(t{+}1)-\psi_w(1)$ for any $w\geq 0$ (compare \eqref{hsthl:falpha}). Note that, if $(P_{\omega,b})_{b\geq 0}$ denotes the contraction semigroup generated by $-T_\omega^+{+}\psi_\omega(1)$ (see Corollary \ref{corspec}), then by Lemma \ref{lmma:Kwasnicki}, \smash{$P_{\omega,b}= \mathit{\Pi}^\ast_\omega e^{-bf_\omega(|\cdot|^2)} \mathit{\Pi_\omega}$} and therefore
\[
P_{\omega,b}\h g(t) = \frac{2}{\pi}\int_0^\infty\hspace{-3pt}\int_0^\infty F_{\omega,\lambda}(t) F_{\omega,\lambda}(s) e^{-bf_\omega(\lambda^2)} g(s)\, ds \, d\lambda
\]
for all $g\in L^1(\mathbb R_+)\cap L^2)(\mathbb R_+)$. It follows that $\int_0^\infty k_{\omega,b}(\cdot,s) g(s) \, ds$ converges in $L^2(\mathbb R_+)$ to $g$ as $b\to 0^+$, where $k_{\omega,b}(t,s)\defeq \frac{2}{\pi}\int_0^\infty F_{\omega,\lambda}(t) F_{\omega,\lambda}(s) e^{-bf_\omega(\lambda^2)} d\lambda$, since $(P_{\omega,b})_{b\geq 0}$ is strongly continuous and $P_{\omega,0}=\mathbb I$. In particular, by changing variables, we obtain for any $\beta>0$ that
\begin{equation} \label{asymhs:eq2}
\beta \int_0^\infty\int_0^\infty \overline{f(t)} \, k_{\omega,b}(\beta t,\beta s) \, g(s) \, ds \h dt \, \xrightarrow{ \, b\to 0^+  } \, (f,g)_{L^2(\mathbb R_+)} 
\end{equation}
for all $f,g\in L^1(\mathbb R_+)\cap L^2(\mathbb R_+)$. And similarly, from $\|P_{w,b} g\|_2 \leq \|g\|_2$ we obtain the uniform bound
\[
\left| \beta \int_0^\infty\int_0^\infty \overline{f(t)} \, k_{\omega,b}(\beta t,\beta s) \, g(s) \, ds \h dt \right| \ \leq \ \|f\|_2\, \|g\|_2 \, ,
\]
which allows the use of dominated convergence below. It follows that
\begin{align*}
& \int_{\mathbb R^d_+} \big(a(\xi){-}1\big)_-^\delta \int_{\mathbb R^d_+} e^{-c|\xi'-\zeta'|^2} e^{-bf_{\mu/|\xi'|}(\zeta_d^2)} \bigg|\int_{\mathbb R^d_+} v_h(\xi,x) \phi(x) \overline{v_h(\zeta,x)} \, dx \bigg|^2 d\zeta \h d\xi \\ 
& = C_{h,d} \int_{\mathbb R^d_+} \big(a(\xi){-}1\big)_-^\delta |\xi'|  \\
& \quad \times  \hspace{-2pt} \int_{\mathbb R^d_+}\int_{\mathbb R^d_+} F_{\mu/|\xi'|,\xi_d}\big(h^{-1}|\xi'|x_d\big) F_{\mu/|\xi'|,\xi_d}\big(h^{-1}|\xi'|y_d\big) \h \phi(x)\phi(y) \\
& \quad \times \hspace{-2pt} \int_{\mathbb R^{d-1}} e^{-c|\xi'-\zeta'|^2-i(\xi'-\zeta')(x'-y')/h} \, |\zeta'| \, k_{\mu/|\zeta'|,b}\Big(\tfrac{|\zeta'|x_d}{h}, \tfrac{|\zeta'| y_d}{h}\Big) \, d\zeta' dy\h dx\h d\xi \\
& = C_{h,d}\int_{\mathbb R^d_+} \big(a(\xi){-}1\big)_-^\delta |\xi'|\int_{\mathbb R^{d-1}}\int_{\mathbb R^{d-1}} \int_{\mathbb R^{d-1}}e^{-c|\xi'-\zeta'|^2 -i(\xi'-\zeta')(x'-y')/h}  \\
& \quad \times  \hspace{-2pt} \int_0^\infty \hspace{-3pt}\int_0^\infty |\zeta'|\, g_{\xi,x'}(x_d)\, k_{\mu/|\zeta'|,b}\Big(\tfrac{|\zeta'|x_d}{h}, \tfrac{|\zeta'| y_d}{h}\Big) g_{\xi,y'}(y_d)\, dy_d \h dx_d \h d\zeta' dy'\h dx'\h d\xi ,
\end{align*}
where $C_{h,d}=4 (2\pi)^{-2d{+}1}  h^{-2d}$ and \smash{$g_{\xi,x'}(x_d):= F_{\mu/|\xi'|,\xi_d}\big(h^{-1}|\xi'| x_d \big)\h \phi(x)$}. By \eqref{asymhs:eq1}, \eqref{asymhs:eq2} and dominated convergence, we obtain
\begin{align*}
& \int_{\mathbb R^d_+} \big(a(\xi){-}1\big)_-^\delta \int_{\mathbb R^d_+} |V\phi V^\ast(\xi,\zeta)|^2 d\zeta \h d\xi \ = \ \frac{h^{-d}}{(2\pi)^{d{-}1}} \frac{2}{\pi} \lim_{c\to 0^+} I_c \, ,
\end{align*}
where 
\begin{align*}
I_c \, &\defeq \ \frac{1}{(2\pi h)^{d-1}} \int_{\mathbb R^d_+} \big(a(\xi){-}1\big)_-^\delta |\xi'| \\
& \, \qquad \times \int_{\mathbb R^{d-1}}\int_{\mathbb R^{d-1}} \int_{\mathbb R^{d-1}} \h e^{-i(\xi'-\zeta')(x'-y')/h} e^{-c|\xi'-\zeta'|^2} d\zeta' \\
& \ \qquad \times \int_0^\infty  F_{\mu/|\xi'|,\xi_d}\big(h^{-1}|\xi'|x_d \big)^2 \h \phi(x',x_d) \h \phi(y',x_d) \, dx_d \h dy'\h dx'\h d\xi \,  \\
& = \, \int_{\mathbb R^d_+} \phi(x) \, \big(\beta_c \ast \phi(\cdot,x_d)\big)(x')\\
& \, \qquad \times \int_{\mathbb R^d_+} \big(a(\xi){-}1\big)_-^\delta |\xi'|\, F_{\mu/|\xi'|,\xi_d}\big(h^{-1}|\xi'|x_d \big)^2 d\xi\, dx,
\end{align*} 
since $(2\pi h)^{-(d-1)}(\mathcal Fe^{-c|\cdot|^2})(\frac{x'{-}y'}{2\pi h}) = \beta_c(x'{-}y')$. Here, \smash{$\beta_c\defeq (\pi/c)^{d/2} e^{-\pi^2|\cdot|^2/c}$} forms an approximate identity in $\mathbb R^{d-1}$ (or \textit{nascent delta function}), in particular \smash{$\int_{\mathbb R^{d-1}} f\h \beta_c$} is uniformly bounded in ${c>0}$, and \smash{$\lim_{c\to 0^+} \int_{\mathbb R^{d-1}} f \h \beta_c = f(0)$} for any \smash{$f\in C(\mathbb R^{d-1})\cap L^\infty(\mathbb R^{d-1})$}. Thus,
\begin{equation}\label{asymhs:eq5050}
\lim_{c\to 0^+ } I_c \, = \, \int_{\mathbb R^d_+} \int_{\mathbb R^d_+} \big(a(\xi){-}1\big)_-^\delta |\xi'|\, F_{\mu/|\xi'|,\xi_d}\big(h^{-1}|\xi'|x_d \big)^2 d\xi  \, \phi(x)^2 \, dx \, .
\end{equation}
In particular, since $\xi\mapsto |\xi'|(a(\xi){-}1)^\delta_-$ belongs to $L^1(\mathbb R^d_+)$ for any $0\leq \delta \leq 1$ and $d\geq 2$, 
\begin{equation}
\int_{\mathbb R^d_+} \big(a(\xi){-}1\big)_-^\delta \int_{\mathbb R^d_+} |V\phi V^\ast(\xi,\zeta)|^2 d\zeta \h d\xi  \, < \, \infty\, , \label{aminusoneHS}
\end{equation}
and thus $(a{-}1)_-^0V\phi V^\ast$ and $(a{-}1)_- V\phi V^\ast$ are Hilbert-Schmidt operators. It follows that
\begin{align*}
\tr\, \phi\h \big(hA^+_{\mu/h}{-}1)_- \phi \, & = \, \tr \, V\phi V^\ast (a{-}1)_-V\phi V^\ast \\
& = \, \int_{\mathbb R^d_+} \big(a(\xi){-}1\big)_- \int_{\mathbb R^d_+} |V\phi V^\ast(\xi,\zeta)|^2 d\zeta \h d\xi \, .
\end{align*}
Hence, \eqref{asymhs:eq3} follows from \eqref{asymhs:eq5050} with $\delta=1$.
\end{proof}

\vspace{10pt}
 
For the upper bound in Proposition \ref{prop:asymphalf} we use

\begin{lem} \label{lmma:asymphalf:upper} Let $\phi\in C_0^1(\mathbb R^d)$ be real-valued  and let \smash{$\rho \defeq \chi (H_{\mu,h}^+)_-^0 \chi$}, where $\chi$ denotes the characteristic function of $\mathrm{supp}(\phi)\cap \mathbb R^d_+$. Then $\rho$ has range in the form domain of $\phi H_{\mu,h}^+\phi$, and for any $\sigma\in (0,\tfrac{1}{2})$,
\begin{align}\nonumber
& \Big|\tr \h \rho \h \phi H_{\mu,h}^+ \phi  \, + \, h^{-d} \Lambda_\mu^{(1)} \int_{\mathbb R^d_+} \phi(x)^2 \, dx \, - \, h^{-d+1} \int_{\mathbb R^d_+} \phi(x)^2 h^{-1} \mathcal K_\mu(h^{-1}x_d)\, dx \Big| \\
& \ \leq \ C_{\sigma} \, (1{+}\mu)^{(d-2\sigma)/2} \,  h^{-d+1+2\sigma} \, . \label{lmma:asymphalf:upper:estimate}
\end{align}
\end{lem}

\begin{proof} Since by definition $\rho f = 0$ in the complement of $\mathbb R^d_+$, similarly as in the proof of Proposition \ref{prop:bulk}, it follows that $\rho f$ belongs to the form domain of $\phi A_{\mu/h}^+ \phi$ for all $f\in L^2(\mathbb R^d)$. Moreover, 
\begin{align} \nonumber
\tr \rho \h \phi H^+_{\mu,h} \phi \, & = \, \tr \rho\h \phi hA_{\mu/h}^+ \phi - \tr \rho\h \phi^2 \\ \nonumber
& = \, \int_{\mathbb R^d_+} \big(a_\mu(\xi){-}1\big)_-^0 \bigg[ \big(\phi^+ \h \overline{v_h(\xi,\cdot)}, hA^+_{\mu/h} \phi^+ \h \overline{v_h(\xi,\cdot)} \big) \\ \label{lmma:asymphalf:upper:1}
& \qquad\qquad\qquad\qquad\qquad - \int_{\mathbb R^d_+}  |v_h(\xi,x)|^2\, \phi(x)^2 \, dx \bigg] \, d\xi \, ,
\end{align}
where $\phi^+:= \chi \phi$. 

For $\varphi \in H_0^{1}(\mathbb R^d_+)$, by the integral representation of the kernel of \smash{$e^{-tA_{\mu/h}^+}$} (see \cite[7.11 -- 7.12]{Lieb2001}), we have
\begin{align*}
(\varphi,A^+_{\mu/h} \varphi)& = \lim_{\delta \to 0^+} (\varphi,A^+_{\mu/h} e^{-\delta A^+_{\mu/h}} \varphi) = - \lim_{\delta\to 0^+} \left.\frac{d}{d\varepsilon}\right|_{\varepsilon=\delta} \big(\varphi, e^{-\varepsilon A^+_{\mu/h}} \varphi\big) \\
& = \, \lim_{\delta,\varepsilon \to 0^+} \int_{\mathbb R^d}\int_{\mathbb R^d} \big|\varphi(x) - \varphi(y)\big|^2\,  \theta^{\delta,\varepsilon}_{\mu/h}(|x{-}y|) \, dx\, dy \, ,
\end{align*}
where for $\nu>0$, $\theta_\nu(t)\defeq \nu^{d+1}\theta(\nu t)$ with $\theta(t)\defeq (2\pi t)^{-(d+1)/2} K_{(d+1)/2}(t)$, and for $\nu,\delta,\varepsilon>0$
\[
\theta_\nu^{\delta,\varepsilon} := \frac{1}{\varepsilon} \Big( (\delta{+}\varepsilon)\h \theta_\nu^{\delta+\varepsilon}- \delta\h \theta_\nu^{\delta}\Big) \ , \quad \theta_\nu^\delta(t) \, \defeq \, \theta_\nu\big((t^2{+}\delta^2)^{1/2}\big) \ , 
\]
in particular $\theta_\nu^0 = \theta_\nu$. Also, we write $\lim_{\delta,\varepsilon\to 0^+}$ to denote the consecutive limits $\lim_{\delta\to 0^+}\lim_{\varepsilon \to 0^+}$, while keeping track of the order of limits. Hence, we have
\begin{align*}
& \Big(\phi^+ \h \overline{v_h(\xi,\cdot)}, A^+_{\mu/h} \phi^+ \h \overline{v_h(\xi,\cdot)} \Big) \\ 
& = \, \lim_{\delta,\varepsilon \to 0^+} \int_{\mathbb R^d_+}\int_{\mathbb R^d_+} \big|\phi(x)\overline{v_h(\xi,h)} - \phi(y)\overline{v_h(\xi,y)} \big|^2 \theta_{\mu/h}^{\delta,\varepsilon}(|x{-}y|) \, dx\, dy \, .
\end{align*}
For fixed $h>0$ and $\xi\in \mathbb R^d_+$, using $v:=v_h(\xi,\cdot)$ as a temporary notation, we write for all $x,y\in \mathbb R^d_+$
\begin{align*}
 & \big|\phi(x)\overline{v_h(\xi,x)} - \phi(y)\overline{v_h(\xi,y)} \big|^2 \\
 & = \, \frac{1}{2}\Big(v(x)\overline{v(y)}(\phi(x){-}\phi(y))^2  + \overline{v(x)}v(y) (\phi(x){-}\phi(y))^2\Big) \\
 &  \quad + \frac{1}{2} \phi(x)^2 \big(2|v(x)|^2 - v(x)\overline{v(y)} - \overline{v(x)}v(y)\big) \\
 &  \quad + \frac{1}{2} \phi(y)^2 \big(2|v(y)|^2 - v(y)\overline{v(x)} - \overline{v(y)}v(x)\big) \, ,
\end{align*}
so that
\begin{equation} \label{lmma:asymphalf:upper:2} 
\big(\phi^+ \h \overline{v_h(\xi,\cdot)}, A^+_{\mu/h} \phi^+ \h \overline{v_h(\xi,\cdot)} \big)  =  \lim_{\delta,\varepsilon \to 0^+}\mathcal R_{\delta,\varepsilon}(\xi) + \lim_{\delta,\varepsilon,\beta\to 0^+}\mathcal I_{\delta,\varepsilon, \beta}(\xi) ,
\end{equation}
where 
\begin{align*}
\mathcal R_{\delta,\varepsilon}(\xi) & \defeq \int_{\mathbb R^d_+}\int_{\mathbb R^d_+} v(x)\overline{v(y)}\, (\phi(x){-}\phi(y))^2 \, \theta^{\delta,\varepsilon}_{\mu/h}(|x{-}y|) \, dx\, dy \, ,\\
\mathcal I_{\delta,\varepsilon,\beta}(\xi) & \defeq \int_{\mathbb R^d_+}\int_{\mathbb R^d_+} \phi(x)^2 \Big(2|v(x)|^2 - v(x)\overline{v(y)} - \overline{v(x)}v(y)\Big)\\
&\qquad \qquad \qquad \qquad \quad \times e_\beta(x)e_\beta(y) \, \theta_{\mu/h}^{\delta,\varepsilon}(|x{-}y|) \, dx\, dy \, ,
\end{align*}
and $e_{\beta}(x):= e^{-\beta_1|x'|^2} e^{-\beta_2 x_d}$, $\beta=(\beta_1,\beta_2)\in \mathbb R_+^2$. As is shown below, the second term in \eqref{lmma:asymphalf:upper:2} combined with the second term in \eqref{lmma:asymphalf:upper:1} yields the two leading terms in the expansion of $\tr \h \rho \h \phi H_{\mu,h}^+ \phi $ stated in the Lemma. Therefore, integrating $\lim_{\delta,\varepsilon\to 0^+}\mathcal R_{\delta,\varepsilon}(\xi)$ in \eqref{lmma:asymphalf:upper:1} results in the remainder, satisfying the estimate \eqref{lmma:asymphalf:upper:estimate}. 

We have
\begin{align*} 
\mathcal I_{\delta,\varepsilon,\beta}(\xi) \, = \, \int_{\mathbb R^d_+}\int_{\mathbb R^d_+}\Big[\frac{1}{2}\big(\mathcal M_\beta + \overline{\mathcal M_\beta}\big)(\xi,x,y) + \mathcal N_\beta(\xi,x,y)\Big] \theta_{\mu/h}^{\delta,\varepsilon}(|x{-}y|) \, dx\, dy \, ,
\end{align*}
where
\begin{align*}
\mathcal M_\beta(\xi,x,y) \, & \defeq \, \Big(\phi(x)^2 v(x) e_\beta(x) - \phi(y)^2 v(y) e_\beta(y)\Big)\Big(\overline{v(x)}e_\beta(x){-}\overline{v(y)}e_\beta(y) \Big) , \\
\mathcal N_\beta(\xi,x,y) \, & \defeq \, \Big(\phi(x)^2|v(x)|^2 e_\beta(x)- \phi(y)^2|v(y)|^2 e_\beta(y)\Big) \Big(e_\beta(y){-}e_\beta(x)\Big) \, .
\end{align*}

First, we show that
\begin{equation} \label{lmma:asymphalf:upper:3}
\lim_{\beta\to 0^+} \int_{\mathbb R^d_+}\int_{\mathbb R^d_+} \mathcal N_\beta(\xi,x,y)\,\theta_{\mu/h}^{\delta,\varepsilon}(|x{-}y|)\, dx \, dy \ = \ 0 \, .
\end{equation}
Pointwise, we have $\lim_{\beta\to 0^+}\mathcal N_\beta(\xi,x,y) = 0$. In order to find a $\beta$-independent integrable upper bound, we separately consider the regions where $|x{-}y|$ is smaller and where $|x{-}y|$ is larger than some $r>0$. For $|x{-}y|>r$, we have 
\begin{align} \nonumber
& \big|\mathcal N_\beta(\xi,x,y)\, \theta^{\delta,\varepsilon}_{\mu/h}(|x{-}y|)\big| \\
& \leq \, C h^{-d} |\xi'| \h\big\|F_{\mu/|\xi'|,\xi_d}\big\|^2_\infty \, \big(\phi(x)^2{+}\phi(y)^2\big)\, |\theta^{\delta,\varepsilon}_{\mu/h}(|x{-}y|)| \, \label{estR:falpha1},
\end{align}
uniformly in $\beta$. This is integrable in the region where $|x{-}y|>r$, because $\phi\in C_0^1(\mathbb R^d)$ and 
\begin{align*}
\int_{|z|>r}\theta^{\delta}_{\mu/h}(|z|) \, dz \, \leq \, C \int_r^\infty t^{-2} dt \, <  \infty \, ,
\end{align*}
where we used that $\theta(t) \leq C \, t^{-(d+1)} e^{-t/2}$ by Lemma \ref{lmma:estbessel} in Appendix \ref{app:bessel}. 

Next, if $|x{-}y|\leq r$, then the condition that $x\in \supp \phi$ or $y\in \supp \phi$ implies that both $x$ and $y$ belong to $B_{R+r}(0)$, where $R>0$ is such that $\supp\phi\subset B_R(0)$. Since
\begin{align*} \big|\phi(x)^2|v(x)|^2-\phi(y)^2|v(y)|^2\big|\, \leq \, Ch^{-1}|\xi'|\h \big\|\nabla \big(\phi^2F_{\mu/|\xi'|,\xi_d}^2 \big)\big\|_\infty\, |x{-}y| 
\end{align*}
and for $\beta_1,\beta_2\leq 1$
\begin{align*} |e_\beta(x)-e_\beta(y)|\ & \leq\ \|\nabla e_\beta\|_\infty \, |x{-}y| \\
& \leq \sup_{x\in\mathbb R^d_+} \big(2\beta_1|x'|e^{-\beta_1|x'|^2}\hspace{-2pt} {+} \h \beta_2 e^{-\beta_2x_d}\big) \, |x{-}y| \leq 3 |x{-}y| \, ,
\end{align*}
it follows for $|x{-}y|\leq r$ that
\begin{equation} \label{estR:falpha2}
\big|\mathcal N_\beta(\xi,x,y)\, \theta^{\delta,\varepsilon}_\mu(|x{-}y|)\big| \,\leq \,  Ch^{-1}|\xi'| \,  \chi_{B_{R+r}(0)}(y)  \, |x{-}y|^2 \, \theta^{\delta,\varepsilon}_{\mu/h}(|x{-}y|)
\end{equation}
uniformly in $\beta_1,\beta_2\leq 1$. The right side is integrable in the region $|x{-}y|<r$, since 
\[ |B_{R+r}(0)| \, \int_{|z|\leq r} |z|^2 \theta^{\delta}_{\mu/h}(|z|) \, dz \ \leq \ C \int_0^r dt \ = \ C r \, < \, \infty \, .\]
Thus, due to \eqref{estR:falpha1} and \eqref{estR:falpha2}, Equation \eqref{lmma:asymphalf:upper:3} follows from dominated convergence. Let $g_h^{\delta,\varepsilon}(x):=\frac{h}{\varepsilon} \, e^{-\delta x/h}(1{-}e^{-\varepsilon x/h})$. By Lemma \ref{lmma:diag},
\begin{align*}
& h \int_{\mathbb R^d_+}\int_{\mathbb R^d_+} \mathcal M_{\beta}(\xi,x,y)\,\theta_{\mu/h}^{\delta,\varepsilon}(|x{-}y|)\, dx \, dy \\
&  = \, \Big(\phi^2 \overline{v_h(\xi,\cdot)} e_\beta, g_h^{\delta,\varepsilon}\big(hA_{\mu/h}^+\big) \overline{v_h(\xi,\cdot)}e_\beta \Big) \\
& = \int_{\mathbb R^d_+}\int_{\mathbb R^d_+}\int_{\mathbb R^d_+} \overline{v_h(\zeta,x)}v_h(\xi,x)\, \phi(x)^2 \\[-5pt]
& \ \, \qquad \qquad \qquad \qquad \qquad \times  e_\beta(x)\, g_h^{\delta,\varepsilon}(a_\mu(\zeta)) \, v_h(\zeta,y) \, \overline{v_h(\xi,y)} \, e_\beta(y) \, dy\, dx\, d\zeta \\[2pt]
& = \frac{|\xi'| h^{-2d}}{(2\pi)^{2(d-1)}} \int_{\mathbb R^d_+}\int_{\mathbb R^d_+}\int_{\mathbb R_+} |\zeta'| \, e^{-ix'(\xi'-\zeta')/h} \, \phi(x)^2\, e_\beta(x) \, g_h^{\delta,\varepsilon}(a_\mu(\zeta)) \\[-2pt]
& \quad \times \tilde F(\xi,\zeta,h^{-1}x_d,h^{-1}y_d)\, e^{-\beta_2 y_d}  \int_{\mathbb R^{d-1}} e^{-iy'(\zeta'-\xi')/h} e^{-\beta_1 |y'|^2} dy' \, dy_d \,  dx \, d\zeta \, , 
\end{align*}
where, for $\xi,\zeta\in \mathbb R^d_+$ and $s,t>0$,
\[
\tilde F(\xi,\zeta,s,t) \defeq  \left(\frac{2}{\pi}\right)^2 f_\xi(s)\, f_\zeta(s) \, f_\xi(t) \, f_\zeta(t) \ , \quad f_\xi(t)\defeq F_{\mu/|\xi'|, \xi_d}(|\xi'| t) \, .
\]
Since $\xi'\mapsto (2\pi h)^{-(d-1)} (\mathcal F^{(d-1)} e^{-\beta_1 |\h \cdot\h |^2})(\xi'/2\pi h)$ defines an approximate identity in $\mathbb R^{d-1}$ with respect to $\beta_1>0$, it follows that
\begin{align*}
&\, h \lim_{\beta_1\to 0^+} \int_{\mathbb R^d_+}\int_{\mathbb R^d_+} \mathcal M_\beta(\xi,x,y)\,\theta_{\mu/h}^{\delta,\varepsilon}(|x{-}y|)\, dx \, dy\\
& = \frac{|\xi'|^2 h^{-d-1}}{(2\pi)^{(d-1)}} \left(\frac{2}{\pi}\right)^2 \int_{\mathbb R^d_+} \phi(x)^2 e^{-\beta_2x_d} \int_{\mathbb R_+} f_\xi(h^{-1} y_d) \, f_\xi(h^{-1} x_d)\, e^{-\beta_2y_d}\\
& \quad \times  \int_{\mathbb R_+} f_{(\xi',\zeta_d)}(h^{-1} y_d)  \, g_h^{\delta,\varepsilon}(a_\mu(\xi',\zeta_d)) \, f_{(\xi',\zeta_d)}(h^{-1} x_d)\, d\zeta_d \, dy_d \, dx \, .
\end{align*}
Hence, by integrating against $\big(a_\mu{-}1\big)_-^0$ (see \eqref{lmma:asymphalf:upper:1}) and changing variables in the $y_d$-integration, we find
\begin{align*}
& C'_{h,d}\, h\int_{\mathbb R^d_+} \big(a_\mu(\xi){-}1\big)_-^0 \lim_{\delta,\varepsilon,\beta\to 0^+} \int_{\mathbb R^d_+}\int_{\mathbb R^d_+} \mathcal M_\beta(\xi,x,y)\,\theta_{\mu/h}^{\delta,\varepsilon}(|x{-}y|)\, dx \, dy \, d\xi \\
& = \lim_{\delta,\varepsilon,\beta_2\to 0^+}\int_{\mathbb R^{d-1}} \int_{\mathbb R^d_+} \int_{\mathbb R_+} \Big(\mathit{\Pi}^\ast_{\mu/|\xi'|}\, |\xi'|\big(a_\mu(\xi',\cdot){-}1\big)_-^0 \, f_{(\xi',\cdot)}(h^{-1}x_d)\Big)(t) \\
& \ \times e^{-h|\xi'|^{-1}\beta_2 t}  \Big(\mathit{\Pi}^\ast_{\mu/|\xi'|}\, g_h^{\delta,\varepsilon}(a_\mu(\xi',\cdot)) \, f_{(\xi',\cdot)} (h^{-1} x_d) \Big)(t) \, dt \ \phi(x)^2\, e^{-\beta_2x_d} \, dx\, d\xi' ,  \\
& = \lim_{\delta,\varepsilon\to 0^+} \int_{\mathbb R^d_+} \int_{\mathbb R^d_+} |\xi'|\big(a_\mu(\xi){-}1\big)_-^0\, g_h^{\delta,\varepsilon}(a_\mu(\xi))\, f_\xi(h^{-1}x_d)^2\, \phi(x)^2 \, dx \, d\xi \, \\
& = C'_{h,d}\int_{\mathbb R^d_+} \int_{\mathbb R^d_+} |\xi'|\big(a_\mu(\xi){-}1\big)_-^0\, a_\mu(\xi)\, |v_h(\xi,x)|^2\, \phi(x)^2 \, dx \, d\xi \, ,
\end{align*}
where $C'_{h,d}\defeq \frac{\pi}{2} (2\pi)^{d-1} h^{d}$, and we are allowed to take limits after the integration in $\xi$ and change the order of integration, since \smash{$\xi \mapsto |\xi'|(a_\mu(\xi){-}1)_-^0$} belongs to $L^1(\mathbb R^d_+)$ whenever $d\geq 2$. 

Considering \eqref{lmma:asymphalf:upper:3}, integrating the second term in \eqref{lmma:asymphalf:upper:2}, and combining the result with the second term in \eqref{lmma:asymphalf:upper:1}, gives
\begin{align*}
& \int_{\mathbb R^d_+} \big(a_\mu(\xi){-}1\big)_-^0 \Big[h \lim_{\delta,\varepsilon,\beta\to 0^+} \mathcal I_{\delta,\varepsilon,\beta}(\xi)- \int_{\mathbb R^d_+} |v_h(\xi,x)|^2\, \phi(x)^2 \, dx \Big] \, d\xi \\
& \ \ = \  \int_{\mathbb R^d_+} \int_{\mathbb R^d_+} |\xi'|\big(a_\mu(\xi){-}1\big)_-\, |v_h(\xi,x)|^2\,  d\xi \ \phi(x)^2 \, dx  \\
& \stackrel{\eqref{asymhs:eq4}}{=} \Lambda^{(1)}_\mu\, h^{-d} \int_{\mathbb R^d_+}\phi(x)^2 dx - h^{-d}\int_{\mathbb R^d_+} \mathcal K_\mu(h^{-1}x_d) \,\phi(x)^2 \h dx \, .
\end{align*}
It remains to prove the bound on the remainder
\begin{align*}
R_{\mu,h}(\phi) & \defeq h \int_{\mathbb R^d_+} \big(a_\mu(\xi){-}1\big)_-^0 \lim_{\delta,\varepsilon \to 0^+} \mathcal R_{\delta,\varepsilon}(\xi) \, d\xi \\
& \, = h \int_{\mathbb R^d_+} \int_{\mathbb R^d_+}\int_{\mathbb R^d_+} \big(a_\mu(\xi){-}1\big)_-^0  v_h(\xi,x)\overline{v_h(\xi,y)}\\
& \ \ \qquad \qquad \qquad \times (\phi(x){-}\phi(y))^2 \, \theta_{\mu/h}(|x{-}y|) \, dx\, dy \, d\xi ,
\end{align*}
where we have used that, by Lemma \ref{lmma:estbessel} and Lemma \ref{app:bessel:deriv} (Appendix \ref{app:bessel}), 
\[
\big|\theta_\nu^{\delta,\varepsilon}(t)\big| \, \leq \, \sup_{\delta\in[0,c]}\left|\frac{d}{d\delta} \big(\delta \theta_\nu^\delta(t)\big)\right| \, \stackrel{}{\leq} \, C_\nu \h t^{-(d+1)} \, ,
\]
and that $\int_{\mathbb R^d_+}\int_{\mathbb R^d_+} (\phi(x){-}\phi(y))^2 |x{-}y|^{-(d+1)} \, dx \h dy \h < \infty$. For $0<\sigma<\frac{1}2$ and $f\in H^{2\sigma}(\mathbb R^{d-1})\cap L^1(\mathbb R^{d-1})$ we have 
\begin{align*}
\frac{|\xi'|^{2\sigma}}{h^{2\sigma}} \int_{\mathbb R^{d-1}} e^{i\xi' x'/h} f(x') \, dx' & = \left(\mathcal F\mathcal F^{-1}|2\pi\cdot|^{2\sigma}\mathcal F f\right)\big(\tfrac{-\xi'}{2\pi h}\big) \\
& 
= \int_{\mathbb R^{d-1}} e^{i\xi'x'/h}\, (-\Delta)^{\sigma}\hspace{-2pt}f\, (x') \, dx' \, ,
\end{align*}
and therefore, by the definition of $v_h(\xi,h)$,
\begin{align*}
R_{\mu,h}(\phi) = h^{1+2\sigma} \int_{\mathbb R^d_+}\int_{\mathbb R^d_+}\int_{\mathbb R^d_+} |\xi'|^{-2\sigma} \big(a_\mu(\xi){-}1\big)_-^0  \, v_h(\xi,x)\,\overline{v_h(\xi,y)} \\ 
\times  \,(-\Delta_{x'})^\sigma \Big(\big(\phi(x){-}\phi(y)\big)^2\theta_{\mu/h}(|x{-}y|)\Big) \, dx\h dy\h d\xi \, .
\end{align*}
Since $|v_h(\xi,x)|\leq C h^{-d/2} |\xi'|^{1/2}$ and
\begin{align*}
& \int_{\mathbb R^d_+}|\xi'|^{1-2\sigma} \big(a_\mu(\xi){-}1\big)_-^0 \, d\xi \,  \leq \,\int_{|\xi'|^2\leq1{+}2\mu}\int_{\xi_d^2\leq (1+2\mu)/|\xi'|^2}|\xi'|^{1-2\sigma} \, d\xi_d \, d\xi' \\
& = \ |\mathbb S^{d-2}| (1{+}2\mu)^{1/2} \int_0^{(1{+}2\mu)^{1/2}} t^{d-2-2\sigma} dt  \, = \,  \frac{|\mathbb S^{d-2}|}{d{-}1{-}2\sigma} (1{+}2\mu)^{(d{-}2\sigma)/2} \, ,
\end{align*}
it follows that
\begin{align*}
& |R_{\mu,h}(\phi)| \\ 
& \leq \, C_\sigma \, \frac{(1{+}\mu)^{(d-2\sigma)/2}}{h^{d-1-2\sigma} } \int_{\mathbb R^d_+}\int_{\mathbb R^d_+} \Big|(-\Delta_{x'})^\sigma \Big(\big(\phi(x){-}\phi(y)\big)^2\theta_{\mu/h}(|x{-}y|)\Big)\Big| \, dx\, dy \, .
\end{align*}
The estimate \eqref{lmma:asymphalf:upper:estimate} now follows from the fact that for any $\sigma\in (0,\frac{1}{2})$ there exists a constant $C_\sigma>0$ such that for all $\nu>0$
\begin{equation}\label{lmma:estLaplacesigma:eq1}
\int_{\mathbb R^d} \int_{\mathbb R^d} \Big|(-\Delta_{x'})^\sigma \Big(\big(\phi(x){-}\phi(y)\big)^2\theta_\nu(|x{-}y|)\Big)\Big| \,dx\h dy \ \leq \  C_\sigma \, .
\end{equation}
The proof of \eqref{lmma:estLaplacesigma:eq1} is a variation of \cite[Lemma B.4]{Frank2016} and is purely technical. It is therefore omitted here (see \cite[Appendix F.2]{Gottwald2016}). \end{proof}

\bigskip 

\begin{proof}[Proof of Proposition \ref{prop:asymphalf}] 
First, by the Variational Principle and Lemma \ref{lmma:hspasymp:01} we have
\begin{align} \label{prop:asymphalf:proof01}
& -\tr\, \big(\phi H^+_{\mu,h}\phi\big)_-  \, \geq  \, -\tr \, \phi\big( H_{\mu,h}^+\big)_- \phi \\
& \ = - h^{-d}\,  \Lambda^{(1)}_\mu \int_{\mathbb R^d_+}\phi(x)^2 dx + h^{-d+1}\int_{\mathbb R^d_+}\phi(x)^2 h^{-1}\mathcal K_\mu(h^{-1}x_d) \, dx \nonumber \, .
\end{align}
Moreover, if $\rho$ is as defined in Lemma \ref{lmma:asymphalf:upper}, then, again by the Variational Principle,
\begin{align} \label{prop:asymphalf:proof02}
& -\tr\, \big(\phi H^+_{\mu,h}\phi\big)_-  \, \leq \, \tr \h \rho \h \phi H_{\mu,h}^+ \phi \\
& = \, - h^{-d} \Lambda_\mu^{(1)} \int_{\mathbb R^d_+} \phi(x)^2 \, dx \, +  \, h^{-d+1} \int_{\mathbb R^d_+} \phi(x)^2 h^{-1} \mathcal K_\mu(h^{-1}x_d)\, dx \, - R_{\mu,h}(\phi) \, , \nonumber
\end{align}
where, by \eqref{lmma:asymphalf:upper:estimate}, for each $\sigma\in (0,\frac{1}{2})$ there exists $C_\sigma>0$ such that
\[
|R_{\mu,h}(\phi)| \, \leq \, C_{\sigma} \, (1{+}\mu)^{(d-2\sigma)/2} \,  h^{-d+1+2\sigma} \, .
\]
Similarly as in \cite[(3.8)]{Frank2013}, recalling that $\Lambda_\mu^{(2)}=\int_0^\infty \mathcal K_\mu(t) \, dt$, we have for any $\delta_1\in (0,1)$
\begin{align*}
&\left|\int_{\mathbb R^d_+} \phi(x)^2 h^{-1} \mathcal K_\mu(h^{-1}x_d)\, dx - \Lambda_\mu^{(2)}\int_{\mathbb R^{d-1}} \phi(x',0)^2 \, dx'\right| \\
& \leq \, \int_0^\infty |\mathcal K_\mu(t)| \left(\int_0^{th}ds\right)^{\delta_1}\left(\int_0^\infty \left|\int_{\mathbb R^{d-1}}\partial_s \phi(x',s)^2 \, dx'\right|^{(1-\delta_1)^{-1}} ds \right)^{1-\delta_1}\\
&\leq \, C_{\delta_1} h^{\delta_1} \int_0^\infty t^{\delta_1} |\mathcal K_\mu(t)| \, dt \ \leq \ C_{\delta_1} (1{+}\mu)^{(d-\delta_1)/2} h^{\delta_1} \, , 
\end{align*}
where the last inequality is due to Lemma \ref{hsthl:lmma:inter1} below. Hence, it follows from \eqref{prop:asymphalf:proof01},
\begin{align}\nonumber
& \tr\, \big(\phi H^+_{\mu,h}\phi\big)_- - h^{-d}\,  \Lambda^{(1)}_\mu \int_{\mathbb R^d_+}\phi(x)^2 dx + h^{-d+1} \Lambda_\mu^{(2)}\int_{\mathbb R^{d-1}} \phi(x',0)^2 \, dx' \\[-2pt] \nonumber
& \leq h^{-d+1}\left|\int_{\mathbb R^d_+}\phi(x)^2 h^{-1}\mathcal K_\mu(h^{-1}x_d) \, dx - \Lambda_\mu^{(2)}\int_{\mathbb R^{d-1}} \phi(x',0)^2 \, dx' \right| \\[2pt]
& \leq \, C_{\delta_1} \, (1{+}\mu)^{(d-\delta_1)/2} \, h^{-d+1+\delta_1} \, , \nonumber
\end{align}
and from \eqref{prop:asymphalf:proof02},
\begin{align}\nonumber
& \tr\, \big(\phi H^+_{\mu,h}\phi\big)_- - h^{-d}\,  \Lambda^{(1)}_\mu \int_{\mathbb R^d_+}\phi(x)^2 dx + h^{-d+1} \Lambda_\mu^{(2)}\int_{\mathbb R^{d-1}} \phi(x',0)^2 \, dx' \\[3pt] \nonumber
& \geq  \, - C_{\delta_2} \, (1{+}\mu)^{(d-\delta_2)/2} h^{-d+1+\delta_2} - C_{\delta_1} \, (1{+}\mu)^{(d-\delta_1)/2} \, h^{-d+1+\delta_1} \, ,
\end{align}
where $\delta_2:=2\sigma$.
\end{proof}

\medskip 

By combining Propositions \ref{prop:straightening} (straightening of the boundary) and \ref{prop:asymphalf} (error in the half-space), we obtain

\begin{cor} \label{prop:asympboundary} There exist constants $c,C>0$ and for all $\delta_1,\delta_2\in (0,1)$ there exist constants $C_{\delta_1},C_{\delta_1,\delta_2}>0$ such that for all real-valued $\phi\in C_0^1(\mathbb R^d)$ satisfying $\|\nabla \phi\|_\infty \leq C l^{-1}$ and supported in a ball of radius $0<l\leq c$ intersecting $\partial\Omega$, 
\begin{align} \label{prop:asympboundary:ineq}
& -C_{\delta_1,\delta_2} \,\bigg( (1{+}\mu)^{(d-\delta_1)/2} \frac{l^{d-1-\delta_1}}{h^{d-1-\delta_1}} + (1{+}\mu)^{(d-\delta_2)/2}\frac{l^{d-1-\delta_2}}{h^{d-1-\delta_2}} \\ 
& \qquad \qquad  +(1{+}\mu)^{d/2} \, w(l)^2\,  \frac{l^{d-1}}{h^{d-1}} + (1{+}\mu)^{d/2} \, w(l)\, \frac{l^{d}}{h^d}\bigg) \nonumber \\
& \leq \,  \tr\, \big(\phi H_{\mu,h}^\Omega \phi\big)_- - h^{-d} \h \Lambda_\mu^{(1)} \int_\Omega \phi(x)^2 dx + h^{-d+1} \h \Lambda_\mu^{(2)} \int_{\partial \Omega} \phi(x)^2 d\sigma(x) \nonumber \\
& \leq \, C_{\delta_1}\left( (1{+}\mu)^{(d-\delta_1)/2} \frac{l^{d-1-\delta_1}}{h^{d-1-\delta_1}} +(1{+}\mu)^{d/2} \, w(l)^2\,  \frac{l^{d-1}}{h^{d-1}} + (1{+}\mu)^{d/2} \, w(l)\, \frac{l^{d}}{h^d}\right) , \nonumber 
\end{align}
where $w$ denotes the modulus of continuity of $\partial \Omega$, see \eqref{modulusofcontinuity}.
\end{cor}
\begin{proof}
From Proposition \ref{prop:straightening} and Proposition \ref{prop:asymphalf}, by rescaling $\phi$, it follows that
\begin{align*}
& \tr\, \big(\phi H_{\mu,h}^\Omega \phi\big)_- - h^{-d} \h \Lambda_\mu^{(1)} \int_\Omega \phi(x)^2 dx + h^{-d+1} \h \Lambda_\mu^{(2)} \int_{\partial \Omega} \phi(x)^2 d\sigma(x) \\
& = \, \tr \, \big(\phi H^+_{\mu,h}\phi\big)_- - h^{-d}\,  \Lambda^{(1)}_\mu \int_{\mathbb R^d_+}\phi(x)^2 dx + h^{-d+1} \Lambda_\mu^{(2)}\int_{\mathbb R^{d-1}} \phi(x',0)^2 \, dx'\\
& \quad \ + \Big(\tr\, \big(\phi H_{\mu,h}^\Omega \phi\big)_- {-}\,  \tr \, \big(\phi H^+_{\mu,h}\phi\big)_- \Big)\\
&\quad \ +  h^{-d+1}\,  \Lambda^{(2)}_\mu \left(\int_{\partial \Omega} \phi(x)^2 d\sigma(x) - \int_{\mathbb R^{d-1}} \phi(x',0)^2 \, dx'  \right)\\[4pt]
& \leq \, C_{\delta_1}\, (1{+}\mu)^{d/2} \big((1{+}\mu)^{-\delta_1/2} (h/l)^{-d+1+\delta_1} + w(l)(h/l)^{-d} +  w(l)^2\, (h/l)^{-d+1} \big) ,
\end{align*}
since, by Lemma \ref{hsthl:lmma:inter1}, $|\Lambda_\mu^{(2)}|\leq C (1{+}\mu)^{d/2}$. Here, we use that rescaling $\phi$ by $\phi_l:=\phi(x/l)$ results in \smash{$\tr(\phi_l H_{\mu,h}^\Omega \phi_l)_- = \tr (\phi H_{\mu,h/l}^\Omega \phi)_-$}, as can be seen by using the integral representation \eqref{quad:rep}. The lower bound follows along the same lines.
\end{proof}

\bigskip

\section{Proof of Theorem \ref{thm:main}} \label{section:proofofthm1}

\begin{proof}
Let $\Omega\subset \mathbb R^d$ be a bounded open domain with $\partial \Omega\in C^1$. First, if $h\geq c$ for some $c>0$, then, for any $\varepsilon >0$, we have
\begin{align*}
& \left|\tr \big(H_{\mu,h}^\Omega\h\big)_- - \Lambda^{(1)}_{\mu} \, |\Omega| \, h^{-d} + \Lambda^{(2)}_{\mu} \, |\partial\Omega|\, h^{-d+1}\right| \\
& \leq \, 2\, \Lambda_\mu^{(1)} |\Omega| \, h^{-d} + \big|\Lambda^{(2)}_{\mu}\big| \, |\partial\Omega|\, h^{-d+1}\, \leq \, C_{\varepsilon} \, (1{+}\mu)^{d/2}\, h^{-d+1+\varepsilon} \, .
\end{align*}
Here, the first inequality follows from $H_{\mu,h}^\Omega = \chi_\Omega H_{\mu,h} \chi_\Omega$ and, by the same argument leading to \eqref{bulk:lemma1:ineq},
\[ 
\tr \, \big(\chi_\Omega H_{\mu,h} \chi_\Omega\big)_- \, \leq \,   \Lambda_\mu^{(1)} \h h^{-d} \int_{\mathbb R^d} \chi_\Omega(x)^2\h dx \, = \, \Lambda_\mu^{(1)} \, |\Omega| \, h^{-d} \, .
\]
In the second inequality we use that \smash{$\Lambda_{\mu}^{(1)}\leq C\, (1{+}\mu)^{d/2}$}, and that, by Lemma \ref{hsthl:lmma:inter1}, \smash{$|\Lambda_\mu^{(2)}|\leq C \, (1{+}\mu)^{d/2}$}. Hence, it remains to prove the claim for small $h$.

For $u\in\mathbb R^d$ let $\phi_u\in C_0^1(\mathbb R^d)$ be as in Section \ref{sec:localization}. By Proposition \ref{prop:loc}, we have
\begin{align} \nonumber
&\tr \h\npb{H_{\mu,h}^\Omega} - \Lambda^{(1)}_\mu |\Omega| \, h^{-d} + \Lambda^{(2)}_\mu |\partial \Omega|\, h^{-d+1} \\ 
& = \int_{\mathbb R^d} \mathcal L_{\mu,h}(\phi_u) \, \frac{du}{l(u)^d} \, + \, \tr \h\npb{H_{\mu,h}^\Omega} - \int_{\mathbb R^d} \tr\h \npb{\phi_u H_{\mu,h}^\Omega \phi_u} \, \frac{du}{l(u)^d} \, , \label{prfmain:eq1}
\end{align}
where
\begin{align*}
\mathcal L_{\mu,h}(\phi_u)\ & \defeq \ \tr\h \npb{\phi_u H_{\mu,h}^\Omega \phi_u} - \Lambda_\mu^{(1)} h^{-d} \int_{\Omega}\phi_u(x)^2dx \\
& \qquad + \Lambda_\mu^{(2)}h^{-d+1} \int_{\partial \Omega}\phi_u(x)^2 d\sigma(x) .
\end{align*}

Note that, if $u\in \mathbb R^d{\setminus}\Omega$ and $\mathrm{supp}(\phi_u)\cap \partial\Omega = \emptyset$, then $\mathcal L_{\mu,h}(\phi_u)=0$. Hence, it suffices to find bounds for $\mathcal L_{\mu,h}(\phi_u)$ when $u$ belongs to the bulk, $u\in U_1\defeq \{u\in\Omega\,|\,B_{l(u)}(u)\cap\partial\Omega = \emptyset\}$, and when $u$ is close to the boundary of $\Omega$, $u\in U_2\defeq \{u\in\mathbb R^d \,  | \, B_{l(u)}(u)\cap \partial \Omega \not = \emptyset\}$. If $u\in U_2$ then it follows from $\delta(u)\leq l(u)$ that $l(u)\leq 3^{-1/2} l_0$. Therefore, by choosing $l_0$ small enough, we are allowed to apply Corollary \ref{prop:asympboundary}. By Proposition \ref{prop:bulk}, in the bulk we have
\[
0 \ \geq \ \int_{U_1} \mathcal L_{\mu,h}(\phi_u) \, \frac{du}{l(u)^d} \ \geq \ -C \, (1{+}\mu)^{(d-1)/2} h^{-d+2} \int_{U_1}l(u)^{-2} \, du \, ,
\]
whereas near the boundary, by Corollary \ref{prop:asympboundary}, for any $\delta,\delta'\in (0,1)$
\begin{align}\nonumber
& -C_{\delta} \, (1{+}\mu)^{d/2} \int_{U_2}\,\left(\frac{l(u)^{-1-\delta}}{h^{d-1-\delta}} +  \frac{w(l(u))^2 l(u)^{-1}}{h^{d-1}} + \frac{w(l(u))}{h^d} \right) \, du \\[2pt] \label{prfmain:ineq1}
& \quad \leq \,  \int_{U_2} \mathcal L_{\mu,h}(\phi_u) \, \frac{du}{l(u)^d} \\[2pt]
& \quad \leq \, C_{\delta'} \, (1{+}\mu)^{d/2} \int_{U_2}\left( \frac{l(u)^{-1-\delta'}}{h^{d-1-\delta'}} + \frac{ w(l(u))^2 l(u)^{-1}}{h^{d-1}} + \frac{w(l(u))}{h^d} \right) \, du \, . \nonumber
\end{align}

By \eqref{intludu}, we have \smash{$\int_{U_1} l(u)^{-2} \, du \, \leq C \, l_0^{-1}$} (since $U_1\subset \Omega^\ast$). Moreover, if $u\in U_2$ then $3^{-1} l_0 < l(u)\leq 3^{-1/2} l_0$. Hence, by the same argument as in the proof of Proposition \ref{prop:loc}, it follows for all $\beta\in \mathbb R$ that $\int_{U_2} l(u)^{\beta} \, du \,\leq C\, l_0^{\beta+1}$. Therefore, by \eqref{prfmain:ineq1}, Proposition \ref{prop:loc}, and \eqref{prfmain:eq1}, for all $h\leq l_0/8$ and $\delta\in (0,1)$
\begin{align} \nonumber
& -C_{\delta} \, (1{+}\mu)^{d/2}\h h^{-d+1} \Big( l_0^{-1} h \, {\mathfrak S}_d(l_0/h) + l_0^{-\delta} h^{\delta} + w(l_0)^2 + w(l_0)\, l_0 \, h^{-1}\Big)  \\[6pt] 
& \leq \,  \tr \h\npb{H_{\mu,h}^\Omega} - \Lambda^{(1)}_\mu |\Omega| \, h^{-d} + \Lambda^{(2)}_\mu |\partial \Omega|\, h^{-d+1} \\[6pt]
& \leq \, C_{\delta} \, (1{+}\mu)^{d/2}\h h^{-d+1} \Big( l_0^{-1} h\, {\mathfrak S}_d(l_0/h) + l_0^{-\delta} h^{\delta} + w(l_0)^2 + w(l_0)\, l_0 \, h^{-1}\Big) \, , \nonumber
\end{align} 
with ${\mathfrak S}_d$ as defined in \eqref{defmathfrakS}.

In the case when $\partial \Omega \in C^{1,\gamma}$, i.e. if $w(t) = C t^{\gamma}$, we choose $l_0$ proportional to $h^{(1+\delta)/(1{+}\delta{+}\gamma)}$ for $d>2$, so that 
\[
 h^{d-1}\left|\tr \h\npb{H_{\mu,h}^\Omega} {-} \,\Lambda^{(1)}_\mu |\Omega| \, h^{-d} {+} \Lambda^{(2)}_\mu |\partial \Omega|\, h^{-d+1}\right| \, \leq \, C_{\delta} \,  (1{+}\mu)^{d/2}\, h^{\delta\gamma/(\gamma{+}1{+}\delta)}.
\]
Since $\varepsilon\defeq\delta\gamma/(\gamma{+}1{+}\delta)$ takes any value in $(0,\gamma/(\gamma{+}2))$ by choosing $\delta \in (0,1)$ appropriately, the error estimate in Theorem \ref{thm:main} follows. In the case of $d=2$, we choose $l_0$ proportional to $h^{2/(\gamma{+}2)}$ and obtain
\[
 h^{d-1}\left|\tr \h\npb{H_{\mu,h}^\Omega} - \Lambda^{(1)}_\mu |\Omega| \, h^{-d} + \Lambda^{(2)}_\mu |\partial \Omega|\, h^{-d+1}\right| \, \leq \, C_{\delta} \,  (1{+}\mu)^{d/2}\, h^\varepsilon
\]
for all $\varepsilon \in (0,\gamma/(\gamma{+}2))$, since $h^{\gamma/(\gamma{+}2)}|\ln(h)|^{1/2} \leq C\, h^{\varepsilon}$ $ \forall \varepsilon \in (0,\gamma/(\gamma{+}2))$.

In the general case of domains with $C^1$ boundaries, let $l_0 = \alpha^{-1} h$, where $\alpha>0$ is such that $8h\leq l_0<\frac{1}{2}$, i.e. $2h<\alpha\leq \frac{1}{8}$. Then, for all $\delta\in(0,1)$, there exists $C>0$ such that
\begin{align*}
r_\mu(h) \, & := \, h^{d-1} (1{+}\mu)^{-d/2} \left|\tr \h\npb{H_{\mu,h}^\Omega} - \Lambda^{(1)}_\mu |\Omega| \, h^{-d} + \Lambda^{(2)}_\mu |\partial \Omega|\, h^{-d+1}\right| \\
& \leq \, C \, \left( \alpha^{\delta} \mathfrak{S}(\alpha^{-1}) + w\left(\frac{h}{\alpha}\right)^2 + \frac{1}{\alpha} \, w\left(\frac{h}{\alpha}\right)  \right) \, \disp{,}
\end{align*}
whenever $0<h<\alpha/2$ and $\mu>0$. 

Let $\varepsilon >0$ and choose $0<\alpha\leq \frac{1}{8}$ such that $\alpha^{\delta}\mathfrak{S}(\alpha^{-1}) < \varepsilon/(2C)$. Then, since $w(t)\to 0$ as $t\to 0^+$, there exists $\delta>0$ such that $h/\alpha<\delta$ implies 
\[
w\big(\tfrac{h}{\alpha}\big)^2 + \tfrac{1}{\alpha} w\big(\tfrac{h}{\alpha}\big) \, < \, \tfrac{\varepsilon}{2C} \, .
\] 
In particular, $r_\mu(h) \, < \, \varepsilon$ for all $h < \min\big\{\alpha/2, \alpha \h \delta\big\}$
and thus $r_\mu(h)\in o(1)$, uniformly in $\mu>0$, as $h\to 0$.
\end{proof}

\bigskip
\section*{Conclusions}

By substituting $h=\lambda^{-1}$, \eqref{thm:asymp:0} is equivalent to the large-$\lambda$ asymptotics of the Riesz mean
\begin{equation} \label{thm:asymp:1}
\sum_{n\in\mathbb N} \big(\lambda_n{-}\lambda\big)_- \ = \ \Lambda_0^{(1)} \, |\Omega| \, \lambda^{d+1} - \Big(\Lambda_0^{(2)}|\partial\Omega|{-}\,C_d \, |\Omega| \, m \Big) \, \lambda^d + \tilde r_m(\lambda) \, ,
\end{equation}

\vspace{-5pt}
\noindent with $\tilde r_m(\lambda) = \, \lambda\, r_m(\lambda^{-1}) \in \mathcal O(\lambda^{d-\varepsilon})$ for any $\varepsilon\in (0,\gamma/(\gamma{+}2))$  when $\partial \Omega\in C^{1,\gamma}$ as $\lambda\to \infty$, and $\tilde r_m(\lambda)\in o(\lambda^d)$ when $\partial \Omega \in C^1$ as $\lambda\to \infty$. Hence, Theorem \ref{thm:asymptotics} is the 
direct generalization of the case $\alpha=1$ in \eqref{frank} for non-zero mass $m>0$. 

Moreover, as is shown in \cite[Lemma A.1]{Frank2013}, from \eqref{thm:asymp:0} we obtain for the $N$-th Ces\`aro mean of the eigenvalues of $A_m^\Omega$,
\vspace{-3pt}
\begin{equation*}
\frac{1}{N}\sum_{n=1}^N \lambda_n \, = \, C_d^{(1)} |\Omega|^{-1/d} N^{1/d} + C_d^{(2)} \Big(\Lambda_0^{(2)}|\partial\Omega| {-}\, C_d\,| \Omega| \, m\Big) |\Omega|^{-1}  + o(1) \, ,
\end{equation*}
as $N\to\infty$, where $C_d^{(1)}=\frac{(d{+}1)^{1+1/d}}{d}\big(\Lambda_0^{(1)}\big)^{-1/d}$ and $C_d^{(2)} = \frac{(d{+}1)^{2d+1}}{d^{2d}} \big(\Lambda_0^{(1)}\big)^{-1}$.

\bigskip
In order to compare with the small-time asymptotics \eqref{asymp:nonlocalheattrace} of the heat trace $Z(t)$ for the eigenvalues of $A_m^\Omega$ by Park and Song \cite{Park2014}, note that the Laplace transform of the map $\lambda\mapsto \sum_n (\lambda_n{-}\lambda)_-$ at $t>0$ is given by $\frac{2}{t^2} Z(t)$. Hence, if $\partial \Omega\in C^{1,\gamma}$, we obtain from \eqref{thm:asymp:1} that for all $\varepsilon\in (\gamma,(\gamma{+}2))$,
\begin{equation}\label{myheattraceresult}
Z(t) = D^{(1)} |\Omega| \, t^{-d}  - \big(D^{(2)} |\partial \Omega| -  D^{(3)} |\Omega|\, m\big) \, t^{-d+1} + \mathcal O(t^{-d+1+\varepsilon})  \, ,
\end{equation}
where $D^{(1)}$, $D^{(2)}$, and $D^{(3)}$ are the constants in \eqref{asymp:nonlocalheattrace} for $\alpha = 1$.
For domains with $C^{1,\gamma}$ boundary, this is a slight improvement upon Park and Song's result, because their remainder is $o(t^{-d+1})$ for Lipschitz domains, and $\mathcal O(t^{-d+2})$ for domains with $C^{1,1}$ boundary.

\begin{rem}[Monotonicity of the subleading term]
Note that the monotonicity in $m$ of the subleading term can be seen already from a purely operator-theoretic argument: If $m>m_0$, then $\psi_m < \psi_{m_0}$ for all $t>0$, and thus $q_m < q_{m_0}$ by \eqref{quad:def}. Hence, by the Variational Principle, we have $\tr(hA_{m}^\Omega{-}1)_- > \tr (hA^\Omega_{m_0}{-}1)_-$, since $(h\psi_m(t)-1)_- > (h\psi_{m_0}(t)-1)_-$ for all $t>0$. Thus, comparing the respective asymptotic expansions in $h$ shows that the coefficients must be monotone in $m$ as well.
\end{rem}

\begin{rem}[Regularity of the boundary] 
Since the contribution to \eqref{theorem1} from a ball intersecting the boundary becomes arbitrarily small when $h\to 0$, it can be shown that our main result extends to Lipschitz domains with boundaries that are $C^1$ except at finitely many points. 

More precisely, if $u_0\in \partial \Omega$, then the support of the corresponding localization function $\phi_{u_0}$ is contained in a ball with radius $l(u_0)\leq \frac{l_0}{2}$, where the localization parameter $l_0$ becomes arbitrarily small when $h\to 0$ (see Sections \ref{sec:localization} and \ref{section:proofofthm1}). Therefore, it can be shown that the contribution from a finite number of points $u_j\in \partial \Omega$, $j=1,\dots, N$, is negligible in the limit $h\to 0$, under the condition that there exist positive constants $R$ and $C$, only depending on the dimension $d$ and $\Omega$, such that
\[
\big| \partial \Omega \cap B_{r}(u_j)\big| \, \leq \, C\, r^{d-1} \, \qquad \forall j=1,\dots, N
\]
for all $r\leq R$. For instance, this condition is satisfied in $d=2$ by any simple polygon.
\end{rem}

\begin{rem}[More general exponents]
Regarding the results of Park and Song \cite{Park2014} and the results of Frank and Geisinger \cite{Frank2013}, it is reasonable to ask whether the approach used in this work can be applied to the operator 
\begin{equation} \label{opalpha}
((-\Delta{+}m^{2/\alpha})^{\alpha/2}{-}m)_D
\end{equation}
with Dirichlet boundary condition on $\Omega$, for arbitrary $\alpha\in (0,2)$. The work \cite{Kwasnicki2011} by Kwa\'snicki, i.e. the explicit diagonalization of the generators of certain L\'evy processes on the half-line, which our method is based on, is also applicable for \eqref{opalpha}. In fact, Kwa\'snicki's diagonalization works for L\'evy processes with L\'evy exponent of the form $f(\xi^2)$, where $f$ is a Bernstein function satisfying $f(0)=0$, and the function $f_{\omega,\alpha}:\mathbb R_+\to \mathbb R_+$, given by
\[
f_{\omega,\alpha}(t) \, \defeq \,  (t{+}1{+}\omega^{2/\alpha})^{\alpha/2}-(1{+}\omega^{2/\alpha})^{\alpha/2} \, ,
\]
is such a Bernstein function for any $\omega>0$ and $\alpha\in (0,2)$ (compare \eqref{hsthl:falpha}). However, our proof of Proposition \ref{prop:straightening} (straightening of the boundary) relies on an integral representation of Modified Bessel functions of the Second Kind (identity \eqref{bessel:intrep1}), which loses the properties we are making use of whenever $\alpha<1$. Other than that, besides a technically more sophisticated analysis of the generalized eigenfunctions of the model operators, there is no reason why the method is not applicable in that case, and of course, it might be possible to prove Proposition \ref{prop:straightening} by other means.
\end{rem}


\bigskip

\appendix

\section{Parallel surfaces of Lipschitz boundaries} \label{app:parallelsurfaces}

For a subset $\Gamma\subset\mathbb R^d$ and $r>0$, the set $\Gamma_r:=\{x\in\mathbb R^d\,|\, \mathrm{dist}(x,\Gamma)< r\}$ is called a \textit{tubular neighbourhood of $\Gamma$}. The boundary $\partial \Gamma_r$ of a tubular neighbourhood is called a \textit{parallel set} (or \textit{surface}) of $\Gamma$.   

\begin{lem} \label{lmma:parallelsurface}
Let $\Omega\subset\mathbb R^d$ be a bounded Lipschitz domain and let $\Gamma:=\partial \Omega$. There exist $\varepsilon >0$ and $C>0$ such that 
\[
\mathcal H^{d-1}(\partial \Gamma_r) \, \leq \, C \qquad \forall r\leq \varepsilon \, ,
\]
where $\mathcal H^{d-1}$ denotes the $(d{-}1)$-dimensional Hausdorff measure.
\end{lem}
\begin{proof}
By \cite[Prop. 5.8]{Kaeenmaeki2013}, for a compact set $\Gamma\subset \mathbb R^d$ there exists a constant $C>0$, such that
\[ 
\mathcal H^{d-1}(\partial \Gamma_r) \ \leq \ C \, r^{d-1} N(\Gamma,r) \qquad \forall r>0,
\]
where $N(\Gamma,r)$ denotes the minimal number of balls of radius $r$ needed to cover $\Gamma$. Clearly, there exists $C>0$ such that $N(\Gamma,r)\leq C \, |\Gamma_r| \, r^{-d}$ for all $r>0$ (see for example \cite[5.6]{Mattila1995}), where $|\Gamma_r|$ denotes the Lebesgue measure of the tubular neighbourhood $\Gamma_r$. The latter is related to the Minkowski $m$-content of $\Gamma$, given by $\mathcal M^m(\Gamma)=\lim_{r \to 0^+}r^{m-d}|\Gamma_r|$ for $0\leq m\leq d$, whenever the limit exists. Since $\Gamma=\partial \Omega$ is a Lipschitz boundary, it is in particular $(d{-}1)$-rectifiable (see e.g. \cite[15.3]{Mattila1995}). Hence, by \cite[3.2.39]{Federer1969}, we have $\mathcal M^{d-1}(\Gamma) = \mathcal H^{d-1}(\Gamma)$, in particular $\lim_{r\to 0^+}|\Gamma_r| r^{-1}$ exists. Therefore there exists $\varepsilon >0$ and $C>0$ such that $|\Gamma_r|\leq C r$ for all $r\leq\varepsilon$, which proves the claim.
\end{proof}

\bigskip
\section{Modified Bessel Functions of the Second Kind} \label{app:bessel}
For $\beta\in \mathbb R$, solutions $s\mapsto K_\beta(s)$ to the Modified Bessel Equation $s^2y''+sy'-(s^2+\beta^2)\,y=0$ are called Modified Bessel Functions of the Second Kind. For $s>0$, as is shown for example in \cite[(9.42)]{Temme1996}, we have
\begin{equation} \label{bessel:intrep00}
K_\beta(s) \, = \, \frac{s^\beta}{2^{\beta+1}} \int_0^\infty e^{-t-s^2/(4t)}\, t^{-\beta-1} dt \, , 
\end{equation}
and by changing variables, see \cite[(9.43)]{Temme1996}, also
\begin{equation} \label{bessel:intrep01}
K_\beta(s) \, = \, \frac{\sqrt{\pi}}{\Gamma\big(\beta{+}\frac{1}{2}\big)} \left( \frac{s}{2} \right)^\beta \int_1^\infty e^{-st} \, (t^2{-}1)^{\nu{-}1/2} dt \, .
\end{equation}

In this work, we are interested in the values $\beta=(n{+}1)/2$ for $n\in\mathbb N$. In that case, these representations yield

\begin{lem}
For any $n\in\mathbb N$, $\nu>0$, and $s>0$, we have
\begin{equation}
K_{(n+1)/2}(\nu s) \, = \, \Big(\frac{s}{\nu}\Big)^{\hspace{-2pt}(n+1)/{2}} \hspace{-5pt} \int_0^\infty  e^{-\nu^2 t -s^2 /(4t)} \, (2t)^{-(n+3)/2} \h dt \,  \label{bessel:intrep2}
\end{equation}
and for any $\alpha \in (0,2]$,
\begin{equation} \label{bessel:intrep2b}
K_{(n+\alpha)/2}\big(\nu^{1/\alpha} s\big) =  \left(\frac{s}{\nu^{1/\alpha}} \right)^{(n+\alpha)/2} \int_0^\infty  e^{-\nu^{2/\alpha} t - s^2/(4t)} (2t)^{-(n+\alpha)/2 - 1} \, dt  \, .
\end{equation}
Moreover, 
\begin{equation}
K_{(n+1)/2} (s) \, = \, \frac{1}{2} \left(\frac{s}{2\pi}\right)^{(n-1)/2} \int_{\mathbb R^n} e^{-s\sqrt{|p|^2 + 1}} \, dp \, .\label{bessel:intrep1}
\end{equation}  
\end{lem}
\begin{proof} The identities \eqref{bessel:intrep2} and \eqref{bessel:intrep2b} follow directly from \eqref{bessel:intrep00} by changing variables. For \eqref{bessel:intrep1}, we note that
\begin{align*}
\int_{\mathbb R^n} e^{-s\sqrt{|p|^2{+}1}} dp \, & = \, |\mathbb S^{n-1}| \int_0^\infty e^{-s\sqrt{r^2+1}} \, r^{n-1} \, dr \\
& = \, \frac{2\pi^{n/2}}{\Gamma(\frac{n}{2})} \, \int_1^\infty e^{-st} (t^2{-}1)^{(n-2)/2}\, t \, dt\\
& = \, \frac{\pi^{n/2}}{\Gamma\big(\frac{n}{2}{+}1\big)} \, \int_1^\infty e^{-st} \frac{d}{dt} (t^2{-}1)^{n/2} \, dt \\
& = \, 2 \left(\frac{2\pi}{s}\right)^{\hspace{-2pt}(n{-}1)/2} \hspace{-3pt} K_{(n+1)/2}(s) \, ,
\end{align*}
where the last equality is due to \eqref{bessel:intrep01}.
\end{proof}

The following lemma is used in the proof of Lemma \ref{lmma:straight1}.  

\begin{lem}\label{lmma:estbessel} For each $n\in\mathbb N_0$ there is a constant $C>0$ such that for all $s>0$
\begin{equation}
K_{(n+1)/2}(s) \, \leq \, C\, s^{-(n+1)/2} e^{-s/2} \, . \label{lmma:estbessel:eq}
\end{equation}
\end{lem}
\begin{proof}
Since $K_{1/2}(s) = \frac{\sqrt{\pi}}{2}\, s^{-1/2}e^{-s} < C s^{-1/2} e^{-s/2}$, the inequality is true for $n=0$. In the case $n\geq 1$, it follows from the integral representation \eqref{bessel:intrep1} and the estimate 
\begin{align*}
\int_{\mathbb R^n} e^{-s\sqrt{|p|^2+1}} dp \ & = \ |\mathbb S^{n-1}| \int_0^\infty \hspace{-6pt}e^{-s\sqrt{t^2+1}} t^{n-1} dt\\
& \leq \ |\mathbb S^{n-1}| \h e^{-s/2} \int_0^\infty \hspace{-6pt}e^{-st/2} t^{n-1} dt \ = \ C \h e^{-s/2}\h s^{-n}\, ,
\end{align*}
that $K_{(n+1)/2}(s) \leq C\, s^{-(n+1)/2} e^{-s/2}$.
\end{proof}

\vspace{6pt}
\begin{lem}\label{lmma:besselineq} For $d\in\mathbb N$ and $s\geq 0$, we have
\begin{equation} 
s \, K_{(d+3)/2}(s)  \leq  2 \h  K_{(d+1)/2}\big(s/\hspace{-1pt}\text{\small$\sqrt{2}$}\big) \, .\label{bessel:ineqrel} \end{equation}
\end{lem}
\begin{proof} We use the integral representation \eqref{bessel:intrep1} with $n=d+2$. For $p\in\mathbb R^{d+2}$, we write $p = (p_d,p_2)$, with $p_d\in\mathbb R^d$ and $p_2\in\mathbb R^2$. Then, $(\sqrt{|p_d|^2{+}1}{-}|p_2|)^2 \geq 0$ implies that
\[
\sqrt{|p|^2+ 1} \ \geq \ \tfrac{1}{\hspace{-2pt}\sqrt{2}} \, \big(\sqrt{|p_d|^2{+}1} + |p_2|\big) \, .
\] 
Hence, by using \eqref{bessel:intrep1} we obtain for $s>0$,
\begin{align*} 
\h K_{(d+3)/2}(s) 
& \, = \, \tfrac{s}{4\pi} \left( \tfrac{s}{2\pi}\right)^{\hspace{-1pt}(d-1)/2} \hspace{-2pt} \int_{\mathbb R^{d+2}} e^{-s\sqrt{|p|^2 + 1}} \, dp  \\
&\leq \tfrac{s}{4\pi} \left( \tfrac{s}{2\pi}\right)^{\hspace{-1pt}(d-1)/2} \hspace{-2pt} \int_{\mathbb R^d} e^{-\frac{s}{\sqrt{2}}\sqrt{|p_d|^2 + 1}} \, dp_d \int_{\mathbb R^2} e^{-\frac{s}{\sqrt{2}}|p_2|}\h dp_2 \\
& \, = \,s\, K_{(d+1)/2}\big(s\hspace{-1pt}/\text{\small$\sqrt{2}$}\big) \int_0^\infty e^{-\frac{s}{\sqrt{2}} r} r \h dr  \, = \,2\h s^{-1} K_{(d+1)/2}\big(s\hspace{-1pt}/\text{\small$\sqrt{2}$}\big)\, .
\end{align*}
Since in the case $s=0$ the inequality \eqref{bessel:ineqrel} is trivially true, this proves the claim.
\end{proof}

\vspace{10pt}
\begin{lem}[Derivative] \label{app:bessel:deriv}For $\beta\in\mathbb R$ and $s>0$, we have 
\begin{equation} \frac{d}{ds} K_{\beta}(s) = \frac{\beta}{s}K_\beta(s) - K_{\beta+1}(s) \, .\end{equation}
\end{lem}
\begin{proof}
This follows immediately from \eqref{bessel:intrep00}, since we are allowed to differentiate under the integral sign, due to
\[
\left|\frac{\partial}{\partial s} \left(e^{-t-s^2/(4t)} \, t^{-\beta-1}\right)\right| = \frac{s}{2}\, e^{-t-s^2/(4t)} \, t^{-\beta-2} \, \leq \, \frac{b}{2} \, e^{-t-a^2/(4t)} t^{-\beta-2}
\]
for all $t>0$ and $s\in [a,b]\subset (0,\infty)$.
\end{proof}

\bigskip

\section{Model operators on the half-line} \label{sec:halfline}

In this section we study the one-dimensional model operators $T_\omega^+$ by applying the results \cite{Kwasnicki2011} by Kwa\'snicki, which provide an explicit spectral representation of the generators of a class of stochastic processes on the half-line. Therefore, in the following we use terminology from the theory of stochastic processes. See \cite[Appendix E]{Gottwald2016} for a concise presentation of the relevant ideas. 

Theorem 1.1 and Theorem 1.3 in \cite{Kwasnicki2011} give a generalized eigenfunction expansion of the generator of a symmetric one-dimensional L\'evy process $X$ killed upon exiting the half-line, with L\'evy exponent of the form $\eta(\xi)= f(\xi^2)$, where $f$ is a so-called \textit{complete Bernstein function} satisfying $\lim_{t\to 0^+}f(t) =0$. Such processes are called \textit{subordinated} to Brownian motion on the real line, which is characterized by the L\'evy exponent $\xi\mapsto \xi^2$. The concept of killing the process when leaving the half-line corresponds to the Dirichlet boundary condition of the associated generator. 

\begin{lem}[\hspace{-0.5pt}{Results from \cite{Kwasnicki2011}}] \label{lmma:Kwasnicki}For a complete Bernstein function $f$ with $f(0+) =0$, let $A$ be the generator in $L^2(\mathbb R_+)$ of the L\'evy process $X$ with L\'evy exponent $\xi\mapsto f(\xi^2)$ killed upon leaving $(0,\infty)$, and let $(P_s)_{s\geq 0}$ denote the contraction semigroup associated to $X$. Then there exists a unitary operator $\mathit{\Pi}$ in $L^2(\mathbb R_+)$ such that $\mathit{\Pi}P_s \mathit{\Pi}^\ast$ is the operator of multiplication by $e^{-sf(\lambda^2))}$ for all $s\geq 0$, and $g\in L^2(\mathbb R_+)$ belongs to $\mathcal D(A)$ if and only if the function $\lambda \mapsto f(\lambda^2)\mathit{\Pi} g (\lambda)$ belongs to $L^2(\mathbb R_+)$, in which case
\begin{equation} \mathit{\Pi} A\h g(\lambda) \, = \, -f(\lambda^2) \mathit{\Pi} g(\lambda) \label{spectralreporig} \end{equation}
for all $\lambda>0$. Moreover, for $\phi\in L^1\cap L^2(\mathbb R_+)$ we have
\begin{equation} \mathit{\Pi} \phi (\lambda) \ = \ \sqrt{\frac{2}{\pi}} \int_0^\infty F_\lambda(t) \h \phi(t)\h dt \, ,\label{defPi}
\end{equation}
where $F_\lambda$ are bounded differentiable functions of the form
\begin{equation}
F_\lambda(t) \, = \,  \sin\big(\lambda t + \vartheta(\lambda)\big) + G_\lambda(t) \label{flambdaexplicitform} \,. 
\end{equation}
Here, $\vartheta$ is given by
\begin{equation} \vartheta (\lambda) \ \defeq \ \frac{1}{\pi} \int_0^\infty \frac{\lambda}{s^2-\lambda^2} \, \ln \frac{(\lambda^2 - s^2)f'(\lambda^2)}{f(\lambda^2)-f(s^2)} \ ds \,  \label{phaseshift}  \, ,
\end{equation}
and $G_\lambda$ is the Laplace transform of a finite measure on $(0,\infty)$, satisfying 
\begin{equation} \label{boundonGlambda}
0 \, \leq \, G_\lambda(t) \, \leq \, \sin\vartheta(\lambda)
\end{equation} 
and
\begin{equation}
\int_0^\infty e^{-t x} \h G_{\lambda} (x)\h dx \ = \ \frac{\lambda \cos \vartheta(\lambda) + t \sin\vartheta(\lambda) }{\lambda^2+ t^2} - \frac{\lambda^2}{\lambda^2{+}t^2} \, \sqrt{\frac{f'(\lambda^2)}{f(\lambda ^2)}} \, \varphi_{\lambda}(t) \, , \label{hsthl:laplacetransfG}
\end{equation}
where
\begin{equation}
\varphi_{\lambda}(t) \ \defeq \ \exp\left( \frac{1}{\pi} \int_0^\infty \frac{t}{t^2{+}s^2}\ln \frac{1-s^2/\lambda^2}{1-f(s^2)/f(\lambda^2)} \ ds \right) \,  
\label{hsthl:philambda}
\end{equation}
for all $t\geq 0$ and $\lambda>0$.
\end{lem}
\begin{proof} 
The main part of the lemma is \cite[Theorem 1.3]{Kwasnicki2011}. Inequality \eqref{boundonGlambda} is proved in \cite[Lemma 4.21]{Kwasnicki2011} and \eqref{hsthl:laplacetransfG} is due to \cite[(4.11)]{Kwasnicki2011} and \cite[Proposition 4.7]{Kwasnicki2011}.
\end{proof}

\vspace{5pt}
\begin{cor}[Spectral representation of $T_\omega^+$] \label{corspec}
For fixed $\omega\geq 0$ and all $\lambda>0$ there exists a real-valued differentiable function $F_{\omega,\lambda}$ on $(0,\infty)$ satisfying $|F_{\omega,\lambda} (x)|\leq 2$ for all $x,\lambda \in (0,\infty)$, such that the operator $\mathit{\Pi}_\omega$ defined by \eqref{defPi} extends to a unitary operator in $L^2(\mathbb R_+)$, and $g\in\mathcal D(T_\omega^+)$ if and only if the function $\lambda \mapsto \psi_\omega(\lambda^2{+}1)\mathit{\Pi}_\omega g (\lambda)$ is in $L^2(\mathbb R_+)$, with $\psi_\omega$ as defined in \eqref{def:psi}. In this case,
\begin{equation} 
\mathit{\Pi_\omega} T_\omega^+ g(\lambda)  \ = \ \psi_\omega(\lambda^2{+}1) \h \mathit{\Pi}_\omega g(\lambda) \,  \label{spectralrepofA} \qquad \textit{for all $\lambda>0$.} \end{equation} 
More precisely, $F_{\omega,\lambda}$ has the form \eqref{flambdaexplicitform}
with the phase shift 
\[
\vartheta_\omega(\lambda) \ = \  \frac{1}{\pi} \int_0^\infty \frac{\lambda}{s^2-\lambda^2} \ \ln \left(\frac{1}{2} + \frac{1}{2} \sqrt{\frac{s^2{+}1{+}\omega^2}{\lambda^2{+}1{+}\omega^2}} \right) \, ds \ = \ \vartheta_0\left(\frac{\lambda}{\sqrt{1{+}\omega^2}}\right) . 
\]
\end{cor}

\begin{proof} We apply Lemma \ref{lmma:Kwasnicki} to the complete Bernstein function $f_\omega$ given by 
\begin{equation}\label{hsthl:falpha}
f_\omega(t)\defeq\psi_\omega(t{+}1)-\psi_\omega(1) \, = \, \sqrt{t{+}1{+}\omega^2}-\sqrt{1{+}\omega^2}\,  \quad \forall t>0\, ,
\end{equation}
and satisfying $\lim_{t\to0^+}f_\omega(t)= f_\omega(0)=0$. If $T$ denotes the $L^2(\mathbb R_+)$ generator of the subordinated L\'evy process with L\'evy symbol $\lambda \to f_\omega(\lambda^2)$ killed upon leaving $(0,\infty)$, then $-T=T_\omega^+  {-}\h \psi_\omega(1)$ (compare \cite[Appendix E.5]{Gottwald2016}). Therefore, \eqref{spectralrepofA} is an immediate consequence of \eqref{spectralreporig}. Moreover, by \eqref{boundonGlambda}, $|F_{\omega,\lambda}(x)| \leq 2$.
\end{proof}

The following lemma provides basic properties of $\vartheta_\omega$ and its first two derivatives.

\begin{lem}[Properties of $\vartheta_\omega$]\label{hspthl:lmma:propvartheta} For each $\omega\geq 0$, the function $\vartheta_\omega$ is monotonically increasing, and twice differentiable on $(0,\infty)$. Moreover,  
\begin{equation}
\frac{d\vartheta_\omega}{d\lambda}(\lambda) \, \leq \, \frac{1}{\pi} \frac{\sqrt{1{+}\omega^2}}{\lambda^2{+}1{+}\omega^2} \ , \ \ 
\left|\frac{d^2\vartheta_\omega}{d\lambda^2}(\lambda) \right| \, \leq \, 
\frac{3}{\pi} \frac{\sqrt{1{+}\omega^2}}{(\lambda^2{+}1{+}\omega^2)^{3/2}}
\label{properties:varthetadiff}
\end{equation}
for all $\lambda>0$, and
\begin{equation}\label{properties:varthetaudiff}
\lim_{\lambda\to 0^+}\vartheta_\omega(\lambda)=0,\quad \lim_{\lambda\to\infty}\vartheta_\omega(\lambda) = \frac{\pi}{8},\quad \lim_{\lambda\to 0^+} \vartheta'_\omega(\lambda) = \frac{1}{\pi\sqrt{1{+}\omega^2}}.
\end{equation}
\end{lem}
\begin{proof} Due to the scaling property $\vartheta_\omega(\lambda) = \vartheta_0(\lambda/\sqrt{1{+}\omega^2})$, we can recover the properties of $\vartheta_\omega$ from those of $\vartheta_0$.

In \cite[Prop 4.17]{Kwasnicki2011} it is proved that, for any complete Bernstein function $f$, the phase shift \eqref{phaseshift} is differentiable, and furthermore that it may be differentiated under the integral sign. Let $l_\omega(s,\lambda)$ denote the logarithm in the definition of $\vartheta_\lambda$. Since
$ \partial_\lambda(\lambda/(s^2{-}\lambda^2)) =(s^2{+}\lambda^2)/(s^2{-}\lambda^2)^2$
is symmetric with respect to an interchange of $\lambda$ and $s$, integrating by parts yields
\begin{align*} \frac{d\vartheta_0}{d\lambda}(\lambda) \ 
& = \ \frac{1}{\pi} \int_0^\infty \left[-\frac{s}{\lambda^2-s^2}\, \frac{\partial}{\partial s} l_0(s,\lambda)  \ + \ \frac{\lambda}{s^2-\lambda^2} \frac{\partial}{\partial\lambda} l_0(s,\lambda) \right] ds \\
& = \ \frac{1}{\pi} \frac{1}{\lambda^2{+}1} \int_0^\infty \Big(\sqrt{s^2{+}1}\Big(\sqrt{\lambda^2{+}1} + \sqrt{s^2{+}1}\Big)\Big)^{-1} ds \, .
\end{align*}
If $t:=s{+}\sqrt{s^2{+}1}$, then $\sqrt{s^2{+}1} = (t^2{+}1)/(2t)$ and $(t^2{+}1)\, dt = 2t\, ds$. Hence, 
\begin{align} \frac{d\vartheta_0}{d\lambda}(\lambda) \ 
& = \ \frac{2}{\pi} \frac{1}{\lambda^2{+}1} \int_{0}^\infty \frac{1}{t^2+2\sqrt{\lambda^2{+}1} \, t +1} \, dt  = \ \frac{1}{\pi} \frac{\tilde l_0(\lambda)}{\lambda(\lambda^2{+}1)}\,  \, , \label{vartheta:derivative}
\end{align}
where for any $\omega\geq 0$, 
\begin{equation} 
\label{hsthl:deflog} \tilde l_\omega(\lambda) \ \defeq \ \ln \frac{\sqrt{\lambda^2{+}1{+}\omega^2}+\sqrt{1{+}\omega^2} + \lambda}{\sqrt{\lambda^2{+}1{+}\omega^2}+\sqrt{1{+}\omega^2} - \lambda} \ = \ \tilde l_0\left( \frac{\lambda}{\sqrt{1{+}\omega^2}}\right) \, .
\end{equation}
Since \eqref{vartheta:derivative} is positive and differentiable in $\lambda>0$, $\vartheta_0$ increases monotonically with $\lambda$ and is twice differentiable. A short calculation shows that 
\begin{equation}\frac{d \tilde l_0}{d\lambda}(\lambda) = \ \frac{1}{\sqrt{\lambda^2{+}1}} \, <  \,  1 . \label{hsthl:inter0} \end{equation}
In particular, since $\tilde l_0(0) = 0$, it follows that $\tilde l_0(\lambda)\leq \lambda$. Thus, \eqref{vartheta:derivative} shows that 
\begin{equation}
\frac{d\vartheta_0}{d\lambda}(\lambda) \, \leq \, \frac{1}{\pi} \frac{1}{\lambda^2{+}1} \, ,
\end{equation}
which implies the first estimate in \eqref{properties:varthetadiff} for any $\omega\geq 0$ by using the scaling property.


By differentiating once more, we find
\begin{align} 
\frac{d^2\vartheta_0}{d\lambda^2}(\lambda) \, & = \, -\frac{1}{\pi \lambda(\lambda^2{+}1)^2} \left( \frac{3\lambda^2 {+} 1}{\lambda} \, \tilde l_0(\lambda) -  \sqrt{\lambda^2{+}1}\right) \disp{.}
\label{prop:vartheta:2ndderivative}
\end{align}
Note that
\[
0 \, \leq \, \frac{d}{d\lambda} \frac{\lambda}{\sqrt{\lambda^2{+}1}} \, = \, \frac{1}{\sqrt{\lambda^2{+}1}}-\frac{\lambda^2}{(\lambda^2{+}1)^{3/2}} \, \leq \, \frac{1}{\sqrt{\lambda^2{+}1}} \, = \, \frac{d \tilde l_0}{d\lambda}(\lambda)\, ,
\]
which (together with $\tilde l_\omega(0)=0$) implies $\tilde l_0(\lambda)  \geq \lambda/{\sqrt{\lambda^2{+}1}}$. In particular, the parenthesis in \eqref{prop:vartheta:2ndderivative} is non-negative for all $\lambda>0$, and therefore
\begin{equation}\left|\frac{d^2\vartheta_0}{d\lambda^2}\right| \ = \ \frac{1}{\pi \lambda(\lambda^2{+}1)^2} \left( \frac{3\lambda^2{+}1}{\lambda} \, \tilde l_0(\lambda) - \sqrt{\lambda^2{+}1}\right) \ \leq \ \frac{3}{\pi} \frac{1}{(\lambda^2{+}1)^{3/2}} \, ,
\end{equation}
where we have used that $\tilde l_0(\lambda) \leq \lambda$ and $\sqrt{\lambda^2{+}1}-1\geq 0$ for all $\lambda\geq 0$. By the scaling property, this completes the proof of \eqref{properties:varthetadiff}. Moreover, it follows from
\[
\frac{1}{\pi} \frac{1}{(\lambda^2{+}1)^{3/2}} \, \leq \, \frac{d\vartheta_0}{d\lambda} \, \leq \, \frac{1}{\pi} \frac{1}{\lambda^2{+}1} \, ,
\]  
that $\lim_{\lambda\to 0^+} \vartheta_\omega'(\lambda) = \frac{1}{\pi} (1{+}\omega^2)^{-1/2}$.

Following \cite[Prop. 4.16]{Kwasnicki2011}, by performing the change of variables $t=s/\lambda$ for $0\,{<}\,s\,{<}\,1$ and $t=\lambda/s$ for $s>1$ in the definition of $\vartheta_0$, we find
\begin{align*} 
\vartheta_0(\lambda) \ 
& = \ \frac{1}{\pi} \int_0^1 \frac{1}{1-t^2} \ \ln \frac{1+\sqrt{\frac{\lambda^2/t^2{+}1}{\lambda^2{+}1}}}{1+\sqrt{\frac{\lambda^2t^2{+}1}{\lambda^2{+}1}}} \ dt \, .
\end{align*} 
By dominated convergence, it follows that $\lim_{\lambda\to 0^+}\vartheta_\omega(\lambda) = 0$, and 
\[ 
\lim_{\lambda \to \infty}  \vartheta_\omega(\lambda) \ = \  \frac{1}{\pi} \int_0^1 \frac{-\ln t}{1-t^2}  \ dt \ = \ \frac{\pi}{8}  \, .
\]
For a proof of the last identity, see for example \cite[Prop. 4.15]{Kwasnicki2011}.
\end{proof}

\vspace{10pt}

Let $G_{\omega,\lambda}$ be the second term in the expression \eqref{flambdaexplicitform} for $F_{\omega,\lambda}$ and let $\varphi_{\omega,\lambda}$ denote the corresponding function \eqref{hsthl:philambda} in the Laplace transform of $G_{\omega,\lambda}$. The following lemma provides properties of $\varphi_{\omega,\lambda}$, which will be needed in the proof of Lemma \ref{hsthl:lmma:inter1} below.

\begin{lem}[Properties of $\varphi_{\omega,\lambda}$]\label{hspthl:lmma:propvarphi} For all $\lambda>0$, the function $\varphi_{\omega,\lambda}$ is differentiable in $t=0$, with
\begin{equation}
\varphi'_{\omega,\lambda}(0) \ = \ \frac{\lambda^2{+}1{+}\omega^2}{1{+}\omega^2}\, \frac{d\vartheta_\omega}{d\lambda}(\lambda) \, ,
\label{hsthl:philambdader}
\end{equation}
and 
\begin{equation}
\lim_{\lambda\to \infty} \varphi'_{\omega,\lambda}(0)  \, = \, 0  \, ,\quad \lim_{\lambda\to 0+} \varphi'_{\omega,\lambda}(0) = \frac{1}{\pi\sqrt{1{+}\omega^2}}\, .\label{hsthl:limphider} 
\end{equation}
\end{lem}
\begin{proof}
If $I_{\lambda,t}(s)$ denotes the integrand in \eqref{hsthl:philambda}, then for any $\varepsilon>0$
\begin{align*}
\frac{1}{\varepsilon} \Big|I_{\lambda,\varepsilon}(s) - I_{\lambda,0}(s)\Big| \ & = \ \frac{1}{\varepsilon^2{+}s^2} \ln \frac{1-s^2/\lambda^2}{1-f_\omega(s^2)/f_\omega(\lambda^2)} \\
& \leq \ \frac{1}{s^2}\ln \frac{1-s^2/\lambda^2}{1-f_\omega(s^2)/f_\omega(\lambda^2)}\eqdef h_\lambda(s)\, .
\end{align*}
We also have
\[
\frac{1-s^2/\lambda^2}{1-f_\omega(s^2)/f_\omega(\lambda^2)} = \frac{f_\omega(\lambda^2)}{\lambda^2} \Big(\hspace{-2pt}\sqrt{s^2 {+}1{+}\omega^2} + \sqrt{\lambda^2 {+}1{+}\omega^2}\,\Big)
\]

\vspace{3pt}
\noindent
and therefore, by l'H\^opital's rule
\[
\lim_{s\to 0+} h_\lambda(s) \, = \, \left[2\Big(\hspace{-2pt}\sqrt{1{+}\omega^2} + \sqrt{\lambda^2 {+}1{+}\omega^2}\,\Big)\sqrt{1{+}\omega^2} \ \right]^{-1} \, .
\]
Hence, $h_\lambda$ is continuous on $[0,\infty)$ and therefore locally integrable near $s=0$. Moreover, since $s\mapsto \ln(s)/s^2$ is integrable on $\left[1,\infty\right)$, $h_\lambda$ is an integrable upper bound for the difference quotient above. Thus, by dominated convergence,
\[
\frac{d\varphi_{\omega,\lambda}}{dt}(0) \, = \,  \frac{1}{\pi} \int_0^\infty \frac{1}{s^2} \, \ln \frac{1-s^2/\lambda^2}{1-f_\omega(s^2)/f_\omega(\lambda^2)} \ ds \, .
\] 
Hence, by monotone convergence, it follows that
\[\lim_{\lambda \to \infty} \frac{d\varphi_{\omega,\lambda}}{dt}(0) \, = \, \frac{1}{\pi} \int_0^\infty \frac{1}{s^2} \lim_{\lambda\to \infty} \ln \frac{1-s^2/\lambda^2}{1-f_\omega(s^2)/f_\omega(\lambda^2)} \h ds\  = \  0 \, ,\]
which proves the first identity in \eqref{hsthl:limphider}. Moreover, integrating by parts yields
\begin{align*}
\frac{d\varphi_{\omega,\lambda}}{dt}(0) \, & = \, \frac{1}{\pi} \int_0^\infty \Big(\sqrt{s^2{+}1{+}\omega^2}\Big(\sqrt{s^2{+}1{+}\omega^2} + \sqrt{\lambda^2{+}1{+}\omega^2}\Big)\Big)^{-1} \, ds \\
&  =\, \frac{\lambda^2{+}1{+}\omega^2}{1{+}\omega^2} \frac{d\vartheta_\omega}{d\lambda}(\lambda) \, ,
\end{align*}
where the last identity follows by comparing with the calculation leading to \eqref{vartheta:derivative}. This shows \eqref{hsthl:philambdader}, and together with \eqref{properties:varthetaudiff} also the second identity in \eqref{hsthl:limphider}.
\end{proof}

\vspace{10pt}
The following result is used in the proof of Proposition \ref{prop:asymphalf} and also provides a bound on the coefficient $\Lambda^{(2)}_\mu = \int_0^\infty \mathcal K_\mu(t) \, dt$ in Theorem \ref{thm:main}, with $\mathcal K_\mu$ as defined in \eqref{defKmu}.

\begin{lem}\label{hsthl:lmma:inter1} For $0 \leq \delta<1$ there exists $C_\delta>0$, such that
\begin{equation}
\int_0^\infty t^\delta |\mathcal K_\mu(t)| \, dt \ \leq \ C_{\delta} \, (1{+}\mu)^{(d-\delta)/2} \, . 
\end{equation}
\end{lem}

\begin{proof}
For any $\nu >0$,
\[
\mathcal J_{\mu,\nu} - \mathcal J_{\mu,\nu}^+ (t)\h  \ = \ \pi^{-1} \int_0^\infty \npb{\psi_{\mu/\nu}(\lambda^2{+}1)-\nu^{-1}} \big(1{-}2F_{\mu/\nu,\lambda}(t)^2\big) \, d\lambda \, ,
\]
where the integrand is non-zero only if $\nu^2(1{+}\lambda^2)\leq 1{+}2\mu$, i.e. $0<\nu\leq \sqrt{1{+}2\mu}$ and 
\begin{equation}
0 \, < \, \lambda \ \leq \left(\frac{1{+}2\mu}{\nu^2} - 1\right)^{1/2} \, \disp{.}
\label{hspthl:rangenu}
\end{equation}
By \eqref{flambdaexplicitform}, we have 
\[
1-2F_{\mu/\nu,\lambda}(t)^2 \ = \ \cos\big(2\beta_{\mu/\nu,t}(\lambda)\big) - 4\sin\big(\beta_{\mu/\nu,t}(\lambda)\big) \, G_{\mu/\nu,\lambda}(t) -2G_{\mu/\nu,\lambda}( t)^2\, ,
\]
where  $\beta_{\omega,t}(\lambda)\defeq \lambda  t + \vartheta_{\omega}(\lambda)$.  Hence, we obtain
\begin{equation}\int_0^\infty t^\delta \left| \mathcal J_{\mu,\nu} -\mathcal J_{\mu,\nu}^+(t)  \right| \h dt \  
\leq \ \pi^{-1}\int_0^\infty t^\delta \Big(|R_1(\nu,t)|+ |R_2(\nu,t)|\Big) \, dt \, , \label{hspthl:step02}
\end{equation}
where
\begin{align*}
R_1(\nu,t) & \defeq \int_0^\infty \Psi_\nu(\lambda) \, \cos\big(2\beta_{\mu/\nu,t}(\lambda)\big) \, d\lambda, \\
R_2(\nu,t) & \defeq \int_0^\infty  \Psi_\nu(\lambda) \, \Big(4\sin\big(\beta_{\mu/\nu,t}(\lambda)\big) \, G_{\mu/\nu,\lambda}( t) +2G_{\mu/\nu,\lambda}( t)^2 \Big)\, d\lambda
\end{align*}
and $\Psi_\nu(\lambda)\defeq \npb{\psi_{\mu/\nu}(\lambda^2{+}1){-}\nu^{-1}}$. Let $0<\delta<1$. We have
\[
\cos\big(2 \beta_{\mu/\nu,t}(\lambda)\big) \ = \ \frac{1}{2t} \left( \frac{d}{d\lambda} \sin\big(2\beta_{\mu/\nu,t}(\lambda)\big) - 2\cos\big(2\beta_{\mu/\nu,t}(\lambda)\big) \frac{d\vartheta_{\mu/\nu}}{d\lambda} \right) \, ,
\]
and integrating by parts in $\lambda$ yields
\begin{align*}
\int_0^1 t^\delta \, |R_1(\nu,t)| \h dt 
& \leq \ \int_0^1 \frac{t^{\delta-1}}{2} \int_0^\Lambda \hspace{-4pt}\left(\left| \frac{d}{d\lambda}\psi_{\mu/\nu}(\lambda^2{+}1) \right| + 2 \nu^{-1} \left|\frac{d\vartheta_{\mu/\nu}}{d\lambda} \right| \right) d\lambda \, dt \, .
\end{align*} 
where $\Lambda = ((1{+}2\mu)/\nu^2{-}1)^{1/2}$. Note that the boundary terms are zero, since by \eqref{properties:varthetaudiff} we have $\lim_{\lambda \to 0^+} \beta_{\mu/\nu,t}(\lambda) = 0$, and $\np{\psi_{\mu/\nu}(\lambda^2{+}1){-}\nu^{-1}}$ vanishes at $\lambda = ((1{+}2\mu)/\nu^2-1)^{1/2}$. We have 
\begin{equation}
\frac{d}{d\lambda} \psi_{\mu/\nu}(\lambda^2{+}1) \, = 
\, \frac{\lambda}{\sqrt{\lambda^2{+}1+(\mu/\nu)^2}} \, ,\label{estimate:etaderivative}
\end{equation}
and thus, 
\begin{align*}
\int_0^\Lambda \left|\frac{d}{d\lambda} \psi_{\mu/\nu}(\lambda^2{+}1)\right| \, d\lambda \, = \,\sqrt{\Lambda^2{+}1{+}(\mu/\nu)^2} - \sqrt{1{+}(\mu/\nu)^2} \, < \, \nu^{-1} \, .
\end{align*}
Moreover, by \eqref{properties:varthetadiff}
\begin{align*} 
\int_0^{\sqrt{1{+}2\mu}/\nu} \left|\frac{d\vartheta_{\mu/\nu}}{d\lambda}\right| \, d\lambda  \, & \leq \,  \frac{1}{\pi}\int_0^{\sqrt{1+2\mu}/\nu} \hspace{-5pt}\frac{\sqrt{1{+}(\mu/\nu)^2}}{\lambda^2{+}1{+}(\mu/\nu)^2}\, d\lambda\\
& \leq \,\frac{1}{\pi}\int_0^\infty \hspace{-5pt}\frac{1}{1{+}x^2} \, dx = \frac{1}{2}  .
\end{align*}
Hence, we obtain
\begin{equation}
\int_0^1 t^\delta \, |R_1(\nu,t)| \h dt \ \leq \ \delta^{-1} \nu^{-1}\label{hspthl:prestep01a} \, .
\end{equation}
In the region where $t\in(1,\infty)$, after two integrations by parts, we find
\begin{align*} 
& \int_1^\infty t^\delta \, |R_1(\nu,t)| \h dt \, \leq \, \int_1^\infty \frac{t^{\delta-2}}{4}\, dt\ \Bigg(\hspace{-2pt}1 + \int_0^{((1{+}2\mu)/\nu^2-1)^{1/2}}\hspace{-4pt} \bigg( \left|\frac{d^2}{d\lambda^2}\psi_{\mu/\nu}(\lambda^2{+}1)\right| \\
& \qquad \quad + \, 3 \left|\frac{d}{d\lambda}\psi_{\mu/\nu}(\lambda^2{+}1)\right| \, \left| \frac{d\vartheta_{\mu/\nu}}{d\lambda}\right| + 2 \nu^{-1} \left|\frac{d\vartheta_{\mu/\nu}}{d\lambda} \right|^2  + \nu^{-1} \left|\frac{d^2\vartheta_{\mu/\nu}}{d\lambda^2}\right|  \Bigg) d\lambda\Bigg), \nonumber
\end{align*}
where we used \eqref{estimate:etaderivative} to bound the non-zero boundary term. We have
\[
\frac{d^2}{d\lambda^2} \h \psi_{\mu/\nu}(\lambda^2{+}1) \ 
= \frac{1+(\mu/\nu)^2}{\sqrt{\lambda^2{+}1{+}(\mu/\nu)^2}^3} \ \leq \ \frac{1}{\sqrt{\lambda^2{+}1{+}(\mu/\nu)^2}} \, ,
\]
and thus, for $\Lambda= {((1{+}2\mu)/\nu^2-1)^{1/2}}$,
\begin{equation}
\int_0^\Lambda \left|\frac{d^2}{d\lambda^2} \psi_{\mu/\nu}(\lambda^2{+}1) \right| \, d\lambda \, \leq \, \int_0^{\Lambda/\sqrt{1{+}{(\mu/\nu)^2}}} \hspace{-10pt}\frac{1}{\sqrt{x^2{+}1}} \ dx \, \leq \,  \frac{1}{\sqrt{\nu^2{+}\mu^2}}\label{estR1:1}\, .
\end{equation}
Next, from \eqref{properties:varthetadiff} and \eqref{estimate:etaderivative}, it follows that
\begin{equation} \int_0^{\sqrt{1+2\mu}/\nu} \left|\frac{d}{d\lambda}\psi_{\mu/\nu}(\lambda^2{+}1)\right|\, \left| \frac{d\vartheta_{\mu/\nu}}{d\lambda}\right| \, d\lambda  \ 
\leq \ \frac{1}{\pi} \int_0^\infty \hspace{-4pt}\frac{1}{x^2{+}1} \ dx \
= \ \frac{1}{2} \, ,
\label{estR1:2}\end{equation}
and 
\begin{align} \nonumber
\nu^{-1} \int_0^{\sqrt{1+2\mu}/\nu} \left|\frac{d\vartheta_{\mu/\nu}}{d\lambda} \right|^2 d\lambda \, & \leq \, \frac{1}{\pi^2\nu} \int_0^{\sqrt{1+2\mu}/\nu} \hspace{-6pt}\frac{1{+}(\mu/\nu)^2}{(1{+}(\mu/\nu)^2{+}\lambda^2)^2} \, d\lambda \\
& \leq \, \frac{1}{2\pi} \frac{1}{\sqrt{\nu^2{+}\mu^2}} \, . \label{estR1:3}
\end{align}
For the last term, by the second estimate in \eqref{properties:varthetadiff}, we obtain
\begin{align} \nonumber
\nu^{-1} \int_0^{\sqrt{1+2\mu}/\nu} \Bigg|\frac{d^2\vartheta_{\mu/\nu}}{d\lambda^2}\Bigg| \, d\lambda \ 
& \leq \ \frac{3}{\pi\nu} \int_0^{\sqrt{1+2\mu}/\nu} \hspace{-6pt}\frac{\sqrt{1{+}(\mu/\nu)^2}}{(\lambda^2{+}1{+}(\mu/\nu)^2)^{3/2}} \, d\lambda \\ 
& \leq  \, \frac{3}{2}\frac{1}{\sqrt{\nu^2{+}\mu^2}} \, .\label{estR1:4}
\end{align}
Combining the estimates \eqref{estR1:1}, \eqref{estR1:2}, \eqref{estR1:3} and \eqref{estR1:4}, 
\[
\int_1^\infty t^\delta \, |R_1(\nu,t)|\h dt \, \leq \, \frac{1}{1{-}\delta} \left(1+\frac{1}{\sqrt{\nu^2{+}\mu^2}}\right)\, .
\]
Together with \eqref{hspthl:prestep01a} this shows that for $0<\delta<1$
\begin{equation} 
\int_0^\infty t^\delta |R_1(\nu,t)|\, dt \ \leq \ C'_\delta \, \Big( 1 + \nu^{-1} \Big) \, , \label{hspthl:step01a}
\end{equation}
where $C'_\delta = 2 \max\{\delta^{-1},(1{-}\delta)^{-1}\}$.

Next, by Lemma \ref{lmma:Kwasnicki}, we have $0\leq G_{\mu/\nu,\lambda}(t)\leq \sin \vartheta_{\mu/\nu}(\lambda)$ for all $t>0$, and therefore 
\begin{align} \nonumber
& \int_0^\infty t^\delta |R_2(\nu,t)| \, dt  \\ 
& \leq \ 6 \int_0^{\sqrt{1+2\mu}/\nu} \hspace{-4pt}\big(\nu^{-1} {-} \psi_{\mu/\nu}(\lambda^2{+}1)\big) \int_0^\infty t^\delta G_{\mu/\nu,\lambda}(t) \, dt \, d\lambda \label{hspthl:R2} \, .
\end{align}
By \eqref{hsthl:laplacetransfG}, for any $\omega>0$
\begin{align*}
\int_0^\infty \hspace{-2pt}G_{\omega,\lambda}(t) \, dt \ & = \ \frac{\cos\vartheta_{\omega}(\lambda)}{\lambda} - \sqrt{\frac{f'_{\omega}(\lambda^2)}{f_{\omega}(\lambda^2)}} \\
& \leq \, \frac{\sqrt{\lambda^2 + \big(\sqrt{\lambda^2{+}1{+}\omega^2}-\sqrt{1{+}\omega^2}\big)^2} - \lambda}{\lambda \sqrt{\lambda^2 + \big(\sqrt{\lambda^2{+}1{+}\omega^2}-\sqrt{1{+}\omega^2}\big)^2}}  \, .
\end{align*}
From here we can perform two different estimates which will be suitable in the cases $\lambda\leq 1$ and $\lambda >1$ respectively. By using
\begin{equation} 
\sqrt{x^2+c^2}-c \, \leq \, \frac{x^2}{2c} \,,  
\quad \sqrt{x^2+c^2}-c \, \leq \, x \, , \quad \forall x,c>0\, ,
\end{equation}
we find
\[
\frac{\sqrt{\lambda^2 + \big(\sqrt{\lambda^2{+}1{+}\omega^2}-\sqrt{1{+}\omega^2}\big)^2} - \lambda}{\lambda \sqrt{\lambda^2 + \big(\sqrt{\lambda^2{+}1{+}\omega^2}-\sqrt{1{+}\omega^2}\big)^2}} \, \leq \ \min\left\{ \frac{\lambda}{8(1{+}\omega^2)}, \frac{1}{\lambda}\right\} \, .
\]
Hence, we obtain for all $\omega,\lambda>0$
\begin{equation} 
\int_0^\infty \hspace{-2pt}G_{\omega,\lambda}(t) \, dt \ \leq \ \min\{\lambda,\lambda^{-1}\} \, . \label{hspthl:int00}
\end{equation}
By differentiating \eqref{hsthl:laplacetransfG}, we also get that for any $\omega,\lambda>0$
\begin{align} 
\int_0^\infty t \h G_{\omega,\lambda} (t)\h dt \, = \, \frac{1}{\lambda} \left(\frac{\tilde l_\omega(\lambda)}{\pi}\sqrt{\frac{f_\omega'(\lambda^2)}{f_\omega(\lambda ^2)}} - \frac{\sin(\vartheta_\omega(\lambda))}{\lambda}\right) \, , \label{hspthl:int01}
\end{align}
where $\tilde l_\omega$ was defined in \eqref{hsthl:deflog}. By Taylor's theorem, there exists $r_\omega\in \mathcal O(1)$ as $\lambda\to0+$, such that
\[
\sin(\vartheta_\omega(\lambda)) \ = \ \cos(\vartheta_\omega(0+)) \, \vartheta'_\omega(0+) \, \lambda +  \lambda^2 \, r_\omega(\lambda) \ = \ \frac{\lambda}{\pi \sqrt{1{+}\omega^2}} +  \lambda^2 \, r_\omega(\lambda) \, ,
\]
since $\vartheta_\omega''(0+) = 0$ by Lemma \ref{hspthl:lmma:propvartheta}. Hence it follows that
\begin{align} 
\int_0^\infty t \h G_{\omega,\lambda} (t)\h dt  \ 
& \leq  \ \frac{1}{\pi\lambda}\left|\tilde l_\omega(\lambda) \sqrt{\frac{f_\omega'(\lambda^2)}{f_\omega(\lambda ^2)}}\, - \frac{1}{\sqrt{1{+}\omega^2}}\right| + |r_\omega(\lambda)|\, . \label{hsthl:inter8} 
\end{align}

Note that, by using the Lagrange form of the remainder, 
for each $\lambda >0$ we can find $\zeta\in (0,\lambda)$ such that 
\begin{align} \nonumber
|r_\omega(\lambda)| \, & = \,  \frac{1}{2} \big|\vartheta'_\omega(\zeta)^2 \sin\vartheta_\omega (\zeta) + \vartheta''_\omega(\zeta)\cos\vartheta_\omega(\zeta) \big| \\ 
& \leq \frac{1}{\pi} \frac{1}{\zeta^2{+}1{+}\omega^2} + \frac{3}{2\pi} \frac{1}{\sqrt{\zeta^2{+}1{+}\omega^2}} \ \leq \ \frac{5}{2\pi} \frac{1}{\sqrt{\zeta^2{+}1{+}\omega^2}} \, , \label{ralpha}
\end{align}
in particular $|r_\omega(\lambda)|<(1{+}\omega^2)^{-1/2}$ for all $\omega,\lambda >0$. 

We proceed by studying \eqref{hsthl:inter8} first for $\lambda<(1{+}\omega^2)^{-1/2}$. Since 
\[\sqrt{\frac{f_\omega(\lambda^2)}{f_\omega'(\lambda ^2)}} \ = 
\ \sqrt{\lambda^2 + \left(\sqrt{\lambda ^2 {+}1{+}\omega^2}{-}\sqrt{1{+}\omega^2}\right)^2} \ \leq  \ \lambda \sqrt{1{+}\lambda^2} \, ,
\]
we have for $0< \lambda \leq (1{+}\omega^2)^{-1/2}$, that
\begin{align}
&\left(\frac{1}{\sqrt{1{+}\omega^2}}{-}\lambda\right)\sqrt{\frac{f_\omega(\lambda^2)}{f_\omega'(\lambda ^2)}} \ \leq \ \frac{1{-}\lambda}{\sqrt{1{+}\omega^2}}\sqrt{\frac{f_\omega(\lambda^2)}{f_\omega'(\lambda ^2)}}  \
 \leq \ \frac{\lambda(1{-}\lambda)\sqrt{1 {+} \lambda^2}}{\sqrt{1{+}\omega^2}}  \,  \label{hsthl:inter9} \, ,
\end{align}
and therefore
\begin{equation}
\tilde l_\omega(\lambda) \ \geq \ \frac{\lambda(1{-}\lambda)\sqrt{1{+}\lambda^2}}{\sqrt{1{+}\omega^2}} \, \geq \, \left(\frac{1}{\sqrt{1{+}\omega^2}}{-}\lambda\right)\sqrt{\frac{f_\omega(\lambda^2)}{f_\omega'(\lambda ^2)}} \,. \label{hsthl:inter20}
\end{equation}
Furthermore, by the proof of Lemma \ref{hspthl:lmma:propvartheta}, $\tilde l_\omega(\lambda) \leq \lambda/\sqrt{1{+}\omega^2}$ for all $\lambda>0$. Thus
\begin{align} \nonumber
\tilde l_\omega(\lambda)\, \sqrt{\frac{f_\omega'(\lambda^2)}{f_\omega(\lambda^2)}} \, & \leq \, \frac{1}{\sqrt{1{+}\omega^2}} \frac{\lambda}{\sqrt{\lambda^2 + \big(\sqrt{\lambda^2{+}1{+}\omega^2}-\sqrt{1{+}\omega^2}\big)^2}} \\ 
& \leq \frac{1}{\sqrt{1{+}\omega^2}} \, , \label{ltildealphatimessqrt}
\end{align}
and therefore, by \eqref{hsthl:inter8},
\[\int_0^\infty t \h G_{\omega,\lambda} (t)\h dt \ \leq  \ \frac{1}{\pi\lambda}\left(\frac{1}{\sqrt{1{+}\omega^2}} - \tilde l_\omega(\lambda) \sqrt{\frac{f_\omega'(\lambda^2)}{f_\omega(\lambda ^2)}} \ \right) + |r_\omega(\lambda)|\,. \]
Together with \eqref{hsthl:inter20} and \eqref{ralpha}, in the case of $0<\lambda\leq (1{+}\omega^2)^{-1/2}$, it follows that
\[
\int_0^\infty t \h G_{\omega,\lambda} (t)\h dt \  \leq \ \frac{1}{\pi} + |r_\omega(\lambda)| \ \leq \ \frac{7}{2\pi} \, .
\]
Next, for $(1{+}\omega^2)^{-1/2}<\lambda\leq 1$, we obtain from \eqref{hsthl:inter8}, \eqref{ralpha} and \eqref{ltildealphatimessqrt} that
\[
\int_0^\infty t \h G_{\omega,\lambda} (t)\h dt \,  
\leq \, \frac{2}{\pi \lambda} \frac{1}{\sqrt{1{+}\omega^2}} + \frac{5}{2\pi}\frac{1}{\sqrt{1{+}\omega^2}} \, < \, \frac{9}{2\pi} \, .
\]
And finally, for $\lambda>1$, by using \eqref{hspthl:int01} and \eqref{ltildealphatimessqrt} we obtain
\[
\int_0^\infty t G_{\omega,\lambda}(t) \h dt \,\leq \, \frac{1}{\lambda} \left(\frac{1}{\pi \sqrt{1{+}\omega^2}} + \frac{1}{\lambda}\right) \, < \, \frac{2}{\lambda} \, .
\] 
Hence, for any $\lambda >0$, we have $\int_0^\infty t \h G_{\omega,\lambda}(t) \, dt <  2 \min\{1,\lambda^{-1}\}$, and thus, by \eqref{hspthl:int00} it follows for all $0\leq \delta < 1$, that
\[
\int_0^\infty t^\delta G_{\mu/\nu,\lambda}(t) \, dt \,  \ \leq \ \int_0^1 G_{\mu/\nu,\lambda}(t) \, dt + \int_1^\infty t \h G_{\mu/\nu,\lambda}(t) \, dt \, < \, 3\min \{1,\lambda^{-1}\} \, .
\]
Thus, by \eqref{hspthl:R2}
\begin{align*}
\int_0^\infty t^\delta |R_2(\nu,t)| \, dt \, & \leq \, 18 \, \nu^{-1} \left(1+ \ln \frac{\sqrt{1{+}2\mu}}{\nu} \right).
\end{align*}
Together with \eqref{hspthl:step01a}, \eqref{hspthl:step02} implies for $0<\delta<1$ that
\begin{align} \ 
& \int_0^\infty t^\delta \left| \mathcal J_{\mu,\nu} - \mathcal J^+_{\mu,\nu}(t)  \right| \h dt \ \leq \ C_\delta \left(1+\nu^{-1} + \nu^{-1} \ln\frac{\sqrt{1{+}2\mu}}{\nu} \right),
\end{align}
and therefore
\begin{align*}
&\int_0^\infty t^\delta|\mathcal K_\mu(t)| \h dt \\
& \leq \, C_{\delta} \, (1{+}2\mu)^{(d-\delta)/2} \left ( \int_0^{1} r^{d-1-\delta} dr + \int_0^1 r^{d-2-\delta} (1{-}\ln r) dr \right) \, .
\end{align*}
Since $\int_0^1 r^{-\beta} dr <\infty$ and $|\int_0^1 r^{-\beta} \ln r \h dr | <\infty$ for any $\beta<1$, it follows that
\[
\int_0^\infty t^\delta |\mathcal K_\mu(t)| \, dt \, \leq \, C_{\delta} \, (1{+}\mu)^{(d-\delta)/2}
\]
for some constant $C_{\delta}>0$ depending only on $\delta\in (0,1)$ and $d\geq 2$. 

In the case $\delta = 0$, the proof follows along the same lines, except for the integration of $|R_1(\nu,t)|$ for small $t$. Here we have
\[
\int_0^\nu |R_1(\nu,t)| \, dt  \, \leq \, \frac{\sqrt{1{+}2\mu}}{\nu}  \, ,
\]
whereas in the region $\nu\leq t\leq \sqrt{1{+}2\mu}$ integration by parts in $\lambda$ yields
\[
\int_\nu^{\sqrt{1{+}2\mu}} |R_1(\nu,t)| \, dt  \, \leq \, \nu^{-1} \int_\nu^1 t^{-1} \, dt  \, = \, \nu^{-1} \ln \frac{\sqrt{1{+}2\mu}}{\nu} \, ,
\]
as in the calculation leading to \eqref{hspthl:prestep01a}. Hence we obtain the same terms as above, and so the claim also follows for $\delta = 0$.
\end{proof}

\bigskip


\vfill

\end{document}